\documentclass[a4paper]{article}
\usepackage[a4paper, total={5.5in, 9in}]{geometry}
\usepackage[alignedleftspaceno]{amsmath}
\usepackage{amsthm,amsfonts,amssymb}
\usepackage[utf8]{inputenc}

\usepackage{tikz}
\usetikzlibrary{calc}
\usetikzlibrary{arrows}
\usetikzlibrary{decorations.markings}

\usepackage{pdfsync}
\usepackage{todonotes}
\usepackage{graphicx}
\usepackage{multirow}
\usepackage{cancel}
\usepackage{subfigure}
\usepackage{euscript}
\usepackage{nicefrac}
\usepackage{enumitem}
\usepackage[bookmarks=true,hidelinks]{hyperref}
\usepackage{cleveref}
\usepackage{bbm}
\usepackage{xcolor}
\usepackage{mathtools}
\usepackage{calrsfs}
\usepackage{esint}

\newtheorem{theorem}{Theorem}
\newtheorem{lemma}[theorem]{Lemma}

\newtheorem{proposition}[theorem]{Proposition}
\newtheorem{corollary}[theorem]{Corollary}

\theoremstyle{definition}
\newtheorem{definition}[theorem]{Definition}
\newtheorem{remark}[theorem]{Remark}
\newtheorem{notation}[theorem]{Notation}
\newtheorem{example}[theorem]{Example}

\numberwithin{theorem}{section}
\numberwithin{equation}{section}

\renewcommand{\geq}{\geqslant}
\renewcommand{\leq}{\leqslant}

\renewcommand{\epsilon}{\varepsilon}
\newcommand{\eps} {\varepsilon}
\renewcommand{\phi}{\varphi}
\newcommand{\R}{\mathbb{R}}

\newcommand{\N}{\mathbb{N}}
\newcommand{\Z}{\mathbb{Z}}
\newcommand{\Lebesgue}{{\EuScript{L}}}	

\newcommand{\Borel}{{\mathcal{B}}}		
\DeclareMathOperator{\sgn}{sgn}

\newcommand{\convhull}[2]{[\![#1,#2]\!]}

\newcommand{\ind}{\mathbbm{1}}
\newcommand{\loc}{{\mathrm{loc}}}
\newcommand{\germ}{\mathcal{G}}
\newcommand{\initialdataset}{L^\infty_{\mathrm{oco}}}		
\newcommand{\constfuncs}{{C_\Dx}}

\newcommand{\Dx}{{\Delta x}}
\newcommand{\Dt}{{\Delta t}}
 \newcommand{\hf}{{\nicefrac{1}{2}}}
 \newcommand{\thf}{{\nicefrac{3}{2}}}

\newcommand{\iphf}{{i+\hf}}
\newcommand{\imhf}{{i-\hf}}

\newcommand{\cell}{{\mathcal{C}}}

\newcommand{\bu}{{\mathbf{u}}}
\newcommand{\bg}{{\mathbf{g}}}
\newcommand{\bv}{{\mathbf{v}}}
\newcommand{\bc}{{\mathbf{c}}}
\newcommand{\bd}{{\mathbf{d}}}

\DeclareMathOperator{\TV}{TV}

\newcommand{\into}{{\mathrm{in}}}
\newcommand{\Nin}{{N_\into}}
\newcommand{\out}{{\mathrm{out}}}
\newcommand{\Nout}{{N_\out}}

\newcommand{\disc}{\mathrm{disc}}
\newcommand{\Ind}{{\mathcal{I}}}			
\newcommand{\Indin}{{\Ind_\into}}
\newcommand{\Indout}{{\Ind_\out}}
\newcommand{\monotone}{{\mathrm{mon}}}
\title{Well-posedness theory for nonlinear scalar conservation laws on networks}
\author{Ulrik Skre Fjordholm, Markus Musch and Nils Henrik Risebro}

\begin{document}
\maketitle
\begin{abstract}
We consider nonlinear scalar conservation laws posed on a network. We establish \(L^1\) stability, and thus uniqueness, for weak solutions satisfying the entropy condition. We apply standard finite volume methods and show stability and convergence to the unique entropy solution, thus establishing existence of a solution in the process. Both our existence and stability/uniqueness theory is centred around families of stationary states for the equation. In one important case -- for monotone fluxes with an upwind difference scheme -- we show that the set of (discrete) stationary solutions is indeed sufficiently large to suit our general theory. We demonstrate the method's properties through several numerical experiments.
\end{abstract}

\section{Introduction}
Partial differential equations  (PDEs) on networks have a large number of applications, including fluid flow in pipelines, traffic flow on a network of roads, blood flow in vessels, data networks, irrigation channels and supply chains.  A treatment of this wide range of applications can be found in the review articles \cite{BressanCanicGaravelloPiccoli2013, Garavello2010} and the references therein. In this paper we will focus on scalar, one-dimensional conservation laws
\begin{equation}\label{eq:cl}
u_t + f(u)_x = 0
\end{equation}
on a network. Here, $u=u(x,t)$ is the unknown \emph{conserved variable} and $f$ is a scalar \emph{flux function} defined either on $\R$ or some subinterval.
We aim to make sense of the conservation law on a directed graph and obtain existence, uniqueness, stability and approximability results.

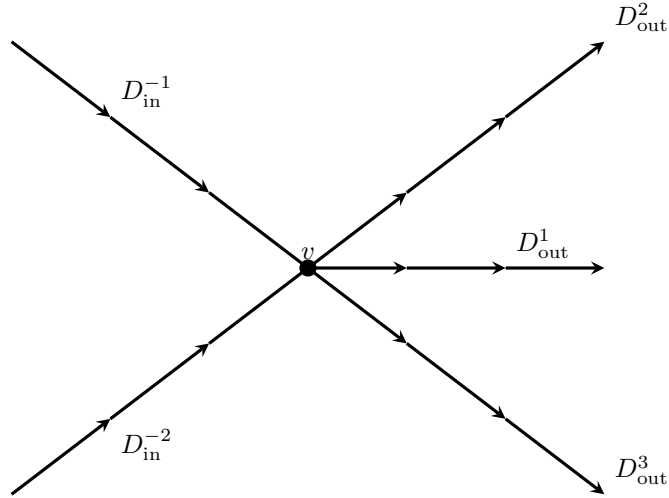
\begin{figure}
  \centering
  \begin{tikzpicture}[>=stealth]
    \coordinate (origin) at (0.0,-3.5); \coordinate (horizontalshift)
    at (0.1,0); \coordinate (verticalshift) at (0,0.1); \coordinate
    (offset1) at (0.5,0); \coordinate (offset2) at (1.5,0);
    \draw[->, very thick] (origin) --
    ($(origin) + 13*(horizontalshift)$); \draw[->, very thick]
    ($(origin) + 13*(horizontalshift)$) --
    ($(origin) -- 26*(horizontalshift)$) node[anchor = south west]
    {$D_\out^1$}; \draw[->, very thick]
    ($(origin) + 26*(horizontalshift)$) --
    ($(origin) -- 39*(horizontalshift)$); \draw[->, very thick]
    ($(origin) + 30*(verticalshift) - 39*(horizontalshift)$) --
    ($(origin) + 20*(verticalshift) - 26*(horizontalshift)$)
    node[anchor = south west] {$D_\into^{-1}$}; \draw[->, very thick]
    (origin) --
    ($(origin) + 13*(horizontalshift) + 10*(verticalshift)$);
    \draw[->, very thick]
    ($(origin) + 13*(horizontalshift) + 10*(verticalshift)$) --
    ($(origin) -- 26*(horizontalshift) + 20*(verticalshift)$);
    \draw[->, very thick]
    ($(origin) + 26*(horizontalshift) + 20*(verticalshift)$) --
    ($(origin) -- 39*(horizontalshift) + 30*(verticalshift)$)
    node[anchor = south west] {$D_\out^2$}; \draw[->, very thick]
    (origin) --
    ($(origin) + 13*(horizontalshift) - 10*(verticalshift)$);
    \draw[->, very thick]
    ($(origin) + 13*(horizontalshift) - 10*(verticalshift)$) --
    ($(origin) -- 26*(horizontalshift) - 20*(verticalshift)$);
    \draw[->, very thick]
    ($(origin) + 26*(horizontalshift) - 20*(verticalshift)$) --
    ($(origin) -- 39*(horizontalshift) - 30*(verticalshift)$)
    node[anchor = south west] {$D_\out^3$}; \draw[->, very thick]
    ($(origin) + 20*(verticalshift) - 26*(horizontalshift)$) --
    ($(origin) + 10*(verticalshift) - 13*(horizontalshift)$);
    \draw[very thick]
    ($(origin) + 10*(verticalshift) - 13*(horizontalshift)$) --
    (origin); \draw[->, very thick]
    ($(origin) - 30*(verticalshift) - 39*(horizontalshift)$) --
    ($(origin) - 20*(verticalshift) - 26*(horizontalshift)$)
    node[anchor = north west] {$D_\into^{-2}$}; \draw[->, very thick]
    ($(origin) - 20*(verticalshift) - 26*(horizontalshift)$) --
    ($(origin) - 10*(verticalshift) - 13*(horizontalshift)$);
    \draw[very thick]
    ($(origin) - 10*(verticalshift) - 13*(horizontalshift)$) --
    (origin); \draw[fill=black] (origin) circle (.7ex) node[anchor =
    south] {$v$};
  \end{tikzpicture}
  \caption{A star shaped network with two ingoing and three outgoing
    edges.}\label{fig:Stargraph}
\end{figure}

Consider a network represented by a connected and directed graph.
We tag the edges of this graph with an index $k$ and impose on each edge a scalar conservation law
\begin{equation}\label{eq:clnetwIntro}
\begin{split}
u_t^k + f^k (u^k)_x = 0&, \qquad x\in D^k,\ t>0 \\
u^k(x,0)=\bar{u}^k(x)&, \qquad x\in D^k
\end{split}
\end{equation}
for some spatial domain $D^k\subset\R$. (Here and in the
remainder, a superscript $k$ will refer to an edge or a vertex.) We
may think of edges as pipes or roads and the vertices as
intersections, with the convention that the direction of travel is in
the positive $x$-direction, as shown in Figure~\ref{fig:Stargraph}.

\begin{example}\label{ex:traffic}
  A common one-dimensional model for traffic flow on a single road is
  described by \eqref{eq:cl} with $f(u)=s_{\max}u(u_{\max}-u)$.  Here
  $u=u(x,t)$ denotes the number of cars per length unit of road
  (density of vehicles), $u_{\max}$ denotes the maximal density of the road, and $s_{\max}$ the maximum speed. Commonly, the length unit is scaled so that $u_{\max}=1$. This model is often referred
  to as the Lighthill--Whitham model \cite{LighthillWhitham1955}.

  A traffic model on a network might incorporate different speed limits
  $s_{\max}^k$ on each road section $k$, or different numbers of
  lanes. In order to model the latter, each unknown $u^k$ can be
  scaled according to the capacity (proportional to the number of lanes) $\alpha^k$
  on road $k$, and the guiding principle to obtain the correct
  scaling should be the conservation of the total mass (i.e., the
  total number of cars) $\sum_k \int_{D^k} u^k(x,t)\,dx$. Thus, if
  $\alpha^k$ denotes the capacity and $v^k(x,t)$ denotes the
  density of vehicles \emph{per lane}, then the
  correct scaling is $u^k=\alpha^kv^k$, and $u^k$ should satisfy
  \eqref{eq:clnetwIntro} with
  \[
    f^k(u^k) = \alpha^k f\big(u^k/\alpha^k\big)
  \]
  and e.g.~$f(v)=v(1-v)$.
\end{example}

In this paper we will be interested in uniqueness and stability for nonlinear scalar conservation laws on a network, as well as in constructing a numerical approximation and proving convergence of the numerical scheme. As opposed to many existing results, where the flux function on each edge is the same \cite{CocliteGaravello2010, HoldenRisebro1995}, we want to allow for a different flux function \( f^k \) on each edge \(D^k\) of the network.
Assuming that  each flux \( f^k \) is continuous and independent of the space variable, the problem can be seen as a PDE with a discontinuous flux, with the points of discontinuity sitting on the vertices of the graph. In fact, if our network would be the trivial network with only one ingoing and one outgoing edge then this would be exactly the case of a conservation law on the real line with a flux function with one discontinuity located at the vertex. Because of the connection to the theory for conservation laws with discontinuous fluxes (see e.g.~\cite{AndreianovKarlsenRisebro2011}), we will borrow several ideas from this theory.

It is well-established that nonlinear hyperbolic conservation laws develop shocks in finite time. Therefore, solutions are always understood in the weak sense.
Unfortunately, weak solutions to nonlinear hyperbolic conservation laws turn out to be non-unique, and additional conditions, usually referred to as entropy conditions, are imposed to select a unique solution. If the flux function is continuous then the theory of entropy solutions is covered by Kruzkhov's theory~\cite{Kruzkov1970}. For conservation laws with discontinuous fluxes the choice of entropy conditions is not obvious, and different physical models might lead to different entropy conditions. Although suitable entropy conditions can yield uniqueness, different entropy conditions are known to yield different solutions; see \cite{AndreianovKarlsenRisebro2011} and references therein. This study of different entropy conditions for conservation laws with discontinuous fluxes culminated in the paper by Andreianov, Karlsen and Risebro~\cite{AndreianovKarlsenRisebro2011}. The authors relate the question of admissibility of a solution to the properties of a set of constant solutions, a so-called germ.
Inspired by the entropy theory of Andreianov, Karlsen and Risebro, we investigate so-called stationary and discrete stationary solutions for our graph problem and thus derive an entropy theory for conservation laws on networks. Although our theory is influenced by the theory in \cite{AndreianovKarlsenRisebro2011}, we have strived to make this paper as self-contained as possible.

In \cite{HoldenRisebro1995} the authors show uniqueness and existence
to the Riemann problem as well as existence of a weak solution of the
Cauchy problem on a network of roads in the case that the flux
functions on each edge are identical. In \cite{CocliteDonadello2020,
  CocliteDiRuvo2019, AndreianovCocliteDonadello2017,
  CocliteGaravello2010} the authors show well-posedness results for
vanishing viscosity solutions. In \cite{GaravelloPiccoli2011} the
authors investigate two entropy type conditions. However, in none of the existing literature can one find a general entropy condition which leads to uniqueness and stability of solutions. In the present work we aim to address this deficiency in the existing theory of conservation laws on networks.

The second important question to address is existence of a solution.
Our approach will be to construct an approximation of the exact entropy solution by constructing a finite volume scheme. We will prove convergence to an entropy solution, thereby also proving existence of a solution.

In this article we focus on monotone fluxes -- that is, each flux $f^k$ is either increasing or decreasing. For non-monotone fluxes only some of the techniques in the present paper are applicable; a traffic flow model, in which fluxes are assumed to be strictly concave, will be treated in a forthcoming work \cite{FjordholmMuschRisebro2020}.

This article is structured as follows: In Section~\ref{Sec:Framework} we define our mathematical framework.
We show uniqueness of entropy solutions to our problem in Section~\ref{Sec:Stability}.
In Section~\ref{Sec:Scheme} we define a finite difference scheme appropriate for our problem, and in Section~\ref{Sec:Convergence} we prove that our numerical scheme converges towards the unique entropy solution.
In Section~\ref{Sec:StationarySol} we show that a class of monotone flux functions fits in our general scheme. Numerical experiments for the monotone case are presented in Section~\ref{Sec:Numerics}.

While the theory outlined in Sections~\ref{Sec:Framework} through~\ref{Sec:Scheme} holds for conservation laws with general flux functions, the convergence theory in Sections~\ref{Sec:Convergence} and~\ref{Sec:Numerics} focuses on monotone flux functions and upwind numerical fluxes.

\section{Entropy solutions on networks}\label{Sec:Framework}
Consider a {network} (or \emph{directed graph}) of vertices and edges; for simplicity we will assume that the network contains a single vertex, along with $N_{\into}\in\N$ edges entering and $N_{\out}\in\N$ edges exiting the vertex (see Figure~\ref{fig:Stargraph}).
(The generalization to general networks will follow analogously, due to the finite speed of propagation of the equations considered.)
We think of the $N_\into$ edges as being to the left of the vertex and the $N_\out$ edges to its right.
The ingoing edges will be labelled \(k\in\Indin\coloneqq\{-\Nin,\dots,-1\}\) and the outgoing edges \(k\in\Indout\coloneqq\{1,\dots,\Nout\}\). We denote \( N \coloneqq \Nin + \Nout \) and we let $\Ind\coloneqq \Indin\cup\Indout$ denote the set of all edge indices.
Placing the vertex at the origin $x=0$, the incoming edges have coordinates $x\in\R_-=(-\infty,0)$, while the outgoing edges have coordinates $x\in\R_+=(0,\infty)$; we will denote the $k$-th edge by
\[
D^k = \begin{cases} \R_- & \text{for } k\in\Indin, \\ \R_+ & \text{for }k\in\Indout.\end{cases}
\]
On each edge $D^k$ we now impose the scalar conservation law \eqref{eq:cl},
resulting in the $N$ distinct equations
\begin{equation}\label{eq:clnetw}
\begin{split}
u^k_t + f^k(u^k)_x = 0&\qquad  \text{for } x\in D^k,\ k \in \Ind.
\end{split}
\end{equation}
The collection of functions \(\bu = (u^k)_{k\in\Ind}\) can be thought of as a function \(\bu\colon\Omega\to\R\), where
\[ \Omega \coloneqq \bigcup_{k\in\Ind} D^k\times\{k\}. \]
On the Borel $\sigma$-algebra \(\Borel(\Omega)\simeq\big\{\prod_{k\in\Ind}A^k\times\{k\}\ :\ A^k\in\Borel(D^k)\big\}\) on $\Omega$ we define the measure $\lambda = \Lebesgue\times \#$, where $\Lebesgue$ is the one-dimensional Lebesgue measure and $\#$ is the counting measure; thus, the integral of \(\bu = (u^k)_{k\in\Ind}\) is
\begin{equation}\label{eq:measuredef}
\int_\Omega \bu\,d\lambda = \sum_{k\in\Ind} \int_{D^k} u^k(x)\,dx.
\end{equation}
The set of $L^\infty$-bounded, real-valued functions on $\Omega$ will be denoted by $L^\infty(\Omega;\lambda)$.
We define the total variation of a function \(\bu \in L^\infty(\Omega;\lambda)\)
as the sum of the variations of its components:
\begin{equation}\label{eq:tvdef}
\TV(\bu) \coloneqq \int_\Omega\Bigl|\frac{d\bu}{dx}\Bigr|\,d\lambda = \sum_{k\in\Ind}\int_{D^k}\Bigl|\frac{du^k}{dx}(x)\Bigr|\,dx.
\end{equation}

\subsection{Weak solutions}
\begin{definition}[Weak Solution]
We say that a function $\bu\in L^\infty\big(\R_+;L^\infty(\Omega;\lambda)\big)$ is a \emph{weak solution} of \eqref{eq:clnetw} with initial data $\bar \bu\in L^\infty(\Omega;\lambda)$ if
\begin{equation}\label{eq:weakform}
\begin{split}
\sum_{k\in\Ind}\int_0^\infty\int_{D^k} u^k\phi^k_t + f^k(u^k)\phi^k_x\,dx\,dt + \sum_{k\in\Ind}\int_{D^k} \bar u^k(x) \phi^k(x,0)\,dx = 0
\end{split}
\end{equation}
for all $\phi^k \in C_c^\infty\big(\overline{D^k}\times[0,\infty)\big)$ satisfying $\phi^k(0,t)\equiv \phi^1(0,t)$ for all $k\in\Ind$.
\end{definition}

Weak solutions automatically satisfy a Rankine--Hugoniot condition at the intersection:

\begin{proposition}[Rankine--Hugoniot condition]\label{prop:RH}
Let $(u^k)_{k\in\Ind}$ be a weak solution of \eqref{eq:clnetw} such that $f^k\circ u^k(\cdot,t)$ has a strong trace at $x=0$ for every $k\in\Ind$ and a.e.~$t>0$. Then
\begin{equation}\label{eq:RH}
\sum_{k\in\Indin}f^k(u^k)(0,t) = \sum_{k\in\Indout}f^k(u^k)(0,t) \qquad\text{for a.e.~} t>0.
\end{equation}
\end{proposition}
\begin{proof}
Define
\begin{equation}\label{hatFct}
\theta_\eps (x) =
\begin{cases}
\frac{1}{\eps} (\eps + x ) & \text{if } x \in [-\eps, 0 ] \\
\frac{1}{\eps} (\eps - x ) & \text{if } x \in [0, \eps ] \\
0                                     & \text{if } |x| > \eps.
\end{cases}
\end{equation}
We define \( \Phi(x,t) \coloneqq \theta_\eps (x) \psi (t) \) where \( \psi \in C_c^\infty ( [ 0, \infty ) ) \).
The partial derivatives of $\Phi$ are
\[
\Phi_x(x, t) =
\begin{cases}
\frac{1}{\eps} \psi (t)  & \text{if } x \in [ -\eps, 0 ] \\
-\frac{1}{\eps} \psi (t) & \text{if } x \in [ 0, \eps ] \\
0                                     & \text{if } |x| > \eps
\end{cases} \qquad \text{and} \qquad
\Phi_t(x,t) = \theta_\eps (x) \psi' (t).
\]
By a density argument, \( \Phi \) qualifies as an admissible test function. Thus, we can insert \( \Phi \) into the weak formulation \eqref{eq:weakform} to get
\begin{align*}
0 ={}& \sum_{k\in\Ind} \int_0^\infty \int_{D^k} u^k \Phi_t^k + f^k (u^k) \Phi_x^k \,dx\,dt + \sum_{k\in\Ind} \int_{D^k} \bar u^k (x) \Phi^k ( x, 0 ) \,dx \\
={}& \sum_{k\in\Ind} \int_0^\infty \int_{D^k} u^k \theta_\eps (x) \psi' (t) \,dt\,dx \\
&+ \frac{1}{\eps} \sum_{k\in\Ind} \int_0^\infty \int_{D^k\cap(-\eps,\eps)} \sgn(k) f^k (u^k) \psi(t) \,dx\,dt \\
&- \frac{1}{\eps} \sum_{k\in\Ind} \int_{D^k\cap(-\eps,\eps)} \bar u^k ( x, 0 ) \psi(t)\,\,dt \\
\to{}&-\sum_{k\in\Ind} \int_0^\infty \sgn(k) f^k (u^k) \psi(t) \,dt
\end{align*}
as $\eps \to 0$, which is equivalent to \eqref{eq:RH}.
\end{proof}
\begin{definition}[Stationary Solution]
A \emph{stationary solution} of \eqref{eq:clnetw} is a weak solution of \eqref{eq:clnetw} which is constant in time and is a \emph{strong} solution on each edge $D^k$.
We see from \eqref{eq:weakform} and \eqref{eq:RH} that the stationary solutions are precisely those satisfying $u^k(x,t)\equiv c^k\in\R$ for $x\in D^k,\ t\geq0$ and $k\in\Ind$, and where $c^k$ satisfy the Rankine--Hugoniot condition
\begin{equation}\label{eqn:RankingeHugoniot} \sum_{k\in\Indin}f^k(c^k) = \sum_{k\in\Indout}f^k(c^k).
\end{equation}
Thus, we can identify each stationary solution with a vector $\bc=(c^k)_{k\in\Ind}\in\R^N$.
\end{definition}
\begin{remark}
Note that if we only required stationary solutions to be \emph{weak} solutions on each edge $D^k$ then they could exhibit arbitrarily many jump discontinuities.
More precisely, if $f$ is not injective then a ``stationary weak solution'' could jump arbitrarily often between values $u^k=c^{k,\pm}$, where $f(c^{k,-})=f(c^{k,+})$.
\end{remark}

\subsection{Entropy conditions}
Next, we formulate conditions that will single out a unique weak solution.
\begin{definition}[Kruzkov entropy pairs]
The \emph{Kruzkov entropy pairs} are the pairs of functions $\eta_c(u)=|u-c|$, $q_c^k(u)=\sgn(u-c)\big(f^k(u)-f^k(c)\big)$ for $c\in\R$.
\end{definition}

The Kruzkov entropy pairs lead to a consistency condition on sets of stationary solutions:
\begin{definition}
A subset $\germ\subset\R^N$ consisting of stationary solutions of \eqref{eq:clnetw} is \emph{mutually consistent} if
\begin{equation}\label{eq:entropyfluxineq}
\sum_{k\in\Indin}q_{c^k}^k\big(\tilde{c}^k\big) \geq \sum_{k\in\Indout}q_{c^k}^k\big(\tilde{c}^k\big)
\end{equation}
for every pair $\bc,\tilde{\bc}\in\germ$, where $q_c^k$ are the Kruzkov entropy flux functions. 
\end{definition}

The set of stationary solutions $\germ$ will determine what class of initial data we can consider:
\begin{definition}\label{Def:Loco}
Let $\germ\subset\R^N$. We let $\initialdataset(\germ)$ be the set of $L^\infty$-bounded data with values lying in the \emph{orthogonal convex hull}\footnote{This nomenclature is somewhat non-standard. The authors found several similar, but non-equivalent, definitions of the orthogonal convex hull in the literature.} of $\germ$,
\begin{equation}\label{eq:definitialdataset}
\begin{split}
\initialdataset(\germ) = \Bigl\{&\bu\in L^\infty(\Omega;\lambda)\ :\ \exists\ \bc,\bd\in\germ \text{ s.t. } [c^k,d^k]\subset\pi^k(\germ)\ \forall\ k\in\Ind \\
&\text{and }
c^k\leq u(x,k)\leq d^k\ \forall\ (x,k)\in\Omega \text{ and } c^k \leq d^l\ \forall\ k,l \in \Ind\Bigr\}
\end{split}
\end{equation}
where $\pi^k$ is the projection $\pi^k(\bc)=c^k$. (Note that we are ignoring the vertex index $k=0$.)
\end{definition}
The condition $[c^k,d^k]\subset\pi^k(\germ)\ \forall\ k\in\Ind$ asserts that for every $\alpha\in[c^k,d^k]$, there is some $\tilde\bc\in\germ$ such that $\tilde{c}^k = \alpha$. This property is needed in the doubling-of-variables argument in the proof of stability and uniqueness.

\begin{example}\label{Ex:Germ}
If $N_\into=N_\out=1$ and $f^k(u)=f(u)=u^2$ then the stationary solutions are all $\bc\in\R^2$ of the form $\bc=(c,c)$ or $\bc=(c,-c)$ for $c\in\R$.
Both $\germ_1=\{(c,c):c\in\R\}$ and $\germ_2=\{(c,-c):c\geq0\}$ (as well as any subset of these) are mutually consistent, as is $\germ_1\cup\germ_2$.
Note that no point of the form $(-c,c)$ for $c>0$ can be added to any of these sets and remain mutually consistent.
Similarly, no set containing both $(-c,c)$ and $(-d,d)$ for distinct $c,d>0$ can be mutually consistent.
The set $\germ_1$ stands out as the smallest closed set which is such that both components span all of $\R$, i.e.~$\pi^k\germ_1=\R$ for $k=-1,1$.
It is readily checked that $\initialdataset(\germ_1) = L^\infty(\Omega;\lambda)$, and that any subset of $\germ_1$ will yield a strictly smaller set of initial data. It is similarly straightforward to check that $\germ_2$ generates a very restrictive set of initial data:
\[
\initialdataset(\germ_2) = \big\{u\in L^\infty(\Omega;\lambda)\ :\ u^{-1}\equiv c,\ u^1\equiv -c \text{ for some } c\geq0\big\}.
\]
Thus, $\germ_1$ is the smallest mutually consistent set of stationary solutions that allow for initial data in all of $L^\infty(\Omega;\lambda)$.

Similar considerations hold when $f^k(u^k)=\alpha^k f(u^k/\alpha^k)$ for $\alpha_{-1}\neq\alpha_1$ (see Example~\ref{ex:traffic}).
\end{example}

\begin{definition}[Entropy Solution]
Let $\germ\subset\R^N$ be a mutually consistent set of stationary solutions of \eqref{eq:clnetw} and let $\bar{\bu}\in L^\infty(\Omega;\lambda)$.
We say that a function $\bu\in L^\infty(\R_+,L^\infty(\Omega;\lambda))$ is an \emph{entropy solution} of \eqref{eq:clnetw} with respect to $\germ$ with initial data $\bar \bu$ if
\begin{equation}\label{eq:entropysoln}
\begin{split}
  \sum_{k\in\Ind}\int_0^\infty\int_{D^k} \eta_{c^k}(u^k)\phi^k_t +
  q_{c^k}^k(u^k)\phi^k_x\,dx\,dt& \\
  + \sum_{k\in\Ind}\int_0^\infty \eta_{c^k}(\bar u^k(x))
  \phi^k(x,0)\,dx &\geq 0
\end{split}
\end{equation}
for every $\bc\in\germ$ and every $0\leq \phi \in C_c^\infty\big(\Omega\times[0,\infty)\big)$ satisfying $\phi^k(0,t)\equiv \phi^1(0,t)$ for all $k\in\Ind$.
\end{definition}
Audusse and Perthame \cite{AudussePerthame2004} considered an entropy condition similar to \eqref{eq:entropysoln}, but in the context of spatially dependent, discontinuous flux functions.

We show first that entropy solutions are invariant in the set $\initialdataset$ from Definition \ref{Def:Loco}.

\begin{lemma}\label{lem:locoinvariant}
Let $\germ\subset\R^N$ be a mutually consistent set of stationary solutions of \eqref{eq:clnetw} and let $\bu\in L^\infty(\R_+,L^\infty(\Omega;\lambda))$ be an entropy solution w.r.t.~$\germ$ with initial data $\bar\bu\in \initialdataset(\germ)$. Then $\bu(t)\in \initialdataset(\germ)$ for a.e.~$t>0$.
\end{lemma}
\begin{proof}
Select $\bc,\bd\in\germ$ such that $\bc\leq \bar{\bu}\leq \bd$ (cf.~\eqref{eq:definitialdataset}). Add inequality \eqref{eq:entropysoln} and equation \eqref{eq:weakform} for both $c^k$ and $u^k$ to get
\begin{align*}
\sum_{k\in\Ind}\int_0^\infty\int_{D^k} \big(c^k-u^k\big)^+\phi^k_t + H\big(c^k-u^k\big)\big(f^k(c^k)-f^k(u^k)\big)\phi^k_x\,dx\,dt& \\
+ \sum_{k\in\Ind}\int_0^\infty \underbrace{\big(c^k-\bar u^k(x)\big)^+}_{=\,0} \phi^k(x,0)\,dx &\geq 0
\end{align*}
(where $\cdot^+=\max(\cdot,0)$ and $H=\sgn^+$ is the Heaviside function). Replacing $\phi$ by a sequence of approximations of the identity function on the set $\Omega\times[0,T]$ yields
\[
\sum_{k\in\Ind}\int_{D^k} \big(c^k-u^k(x,T)\big)^+\,dx \leq 0
\]
for a.e.~$T>0$, whence $\bu(T)\geq\bc$. It follows similarly that $\bu(T)\leq\bd$, and hence, $\bu(T)\in\initialdataset(\germ)$.
\end{proof}

The above lemma enables us to show that entropy solutions have strong traces.

\begin{theorem}
	Let $\germ\subset\R^N$ be a mutually consistent set of stationary solutions of \eqref{eq:clnetw} and let $\bu\in L^\infty(\R_+,L^\infty(\Omega;\lambda))$ be an entropy solution w.r.t.~$\germ$ with initial data $\bar\bu\in \initialdataset(\germ)$. Then the functions \( q^k\circ u^k \) and $f^k\circ u^k$ admit strong traces on \( \{x = 0\}\), for any $k\in\Ind$.
\end{theorem}
\begin{proof}
	It follows from Lemma \ref{lem:locoinvariant} that $u^k$ is a Kruzkhov entropy solution on $D^k$, for any $k\in\Ind$. We can therefore apply \cite[Theorem 1.4]{Panov2007} to obtain the desired conclusion.
\end{proof}

\begin{proposition}\label{prop:decreasingEntropy}
Let $\germ\subset\R^N$ be a set of stationary solutions of \eqref{eq:clnetw}. Let \(\bu\) be an entropy solution of \eqref{eq:clnetw} w.r.t.~$\germ$ such that $q_{c^k}^k\circ u^k(\cdot,t)$ has a strong trace at $x=0$ for every \(k\in\Ind\) and a.e.~$t>0$. Then
\begin{equation}\label{eq:RHentr}
\sum_{k\in\Indin}q_{c^k}^k(u^k)(0,t) \geq \sum_{k\in\Indout}q_{c^k}^k(u^k)(0,t) \qquad \text{for a.e.~}t>0
\end{equation}
for every $\bc\in\germ$.
\end{proposition}
\begin{proof}
We take a positive test function $0 \leq \psi \in C_c^\infty((0, \infty))$.
As in the proof of Proposition~\ref{prop:RH} we define $\Phi(x, t) \coloneqq \theta_\eps (x,t) \psi (t)$ where $\theta_\eps$ is given by \eqref{hatFct}.
Now we insert $\Phi$ as test function into the entropy inequality \eqref{eq:entropysoln} to get
\begin{align*}
  0 \leq{}& \sum_{k\in\Ind} \int_0^\infty \int_{D^k} \eta_{c^k} (u^k) \theta_\eps (x) \psi' (t) \,dx\,dt \\
    &-\frac{1}{\eps} \sum_{k\in\Ind} \biggl(\begin{aligned}[t]&\int_0^\infty \int_{D^k\cap(-\eps,\eps)}  \sgn(k) q_{c^k}^k (u^k) \psi (t) \,dx\,dt \\
    &+ \int_{D^k} \eta_{c^k} ( \bar u^k (x) ) \theta_\eps (x) \psi (t) \,dx \biggr)\end{aligned} \\
    \to{}&- \sum_{k\in\Ind} \int_0^\infty \sgn(k) q_{c^k}^k ( u^k(0,t) ) \psi (t) \,dt
\end{align*}
as $\eps\to0$, which shows the desired inequality.
\end{proof}

\begin{definition}[Maximality]
A mutually consistent set of stationary solutions \( \germ \) is called \emph{maximal} if there is no mutually consistent set of stationary solutions $\widetilde{\germ}$ having $\germ$ as a strict subset.
\end{definition}

\begin{lemma}
Let $\germ$ be a maximal, mutually consistent set of stationary solutions and let \( \bd \in \R^{N_\into + N_\out} \) satisfy the Rankine--Hugoniot condition \eqref{eqn:RankingeHugoniot} as well as
\[
\sum_{k\in\Indin}q_{c^k}^k\big(d^k\big) \geq \sum_{k\in\Indout}q_{c^k}^k\big(d^k\big) \qquad\forall\ \bc\in\germ.
\]
Then $\bd \in \germ$.
\end{lemma}
\begin{proof}
By assumption, $\germ\cup\{\bd\}$ is a mutually consistent set of stationary solutions. But $\germ$ is maximal, whence $\germ\cup\{\bd\}=\germ$.
\end{proof}

\begin{remark}
By Proposition \ref{prop:decreasingEntropy}, we can take the vector \( \bd \) to be the trace of an entropy solution at the vertex. This observation will be important in the proof of the stability result in Section \ref{Sec:Stability},
\end{remark}

\section{Stability and Uniqueness}\label{Sec:Stability}
\begin{theorem}[Entropy Solutions are \(L^1\) stable]\label{thm:entrsolnunique}
Let $\germ\subset\R^N$ be a mutually consistent, maximal set of stationary solutions.
Let $\bu,\,\bv$, be entropy solutions of \eqref{eq:clnetw} w.r.t.~$\germ$ with initial data \(\bar\bu,\, \bar\bv \in \initialdataset(\germ)\cap L^1(\Omega;\lambda)\), and assume each function $u^k,v^k$ has a strong trace at $x=0$.
Then
\[
\sum_{k\in\Ind} \big\|u^k(t)-v^k(t)\big\|_{L^1(D^k)} \leq \sum_{k\in\Ind} \big\|\bar u^k-\bar v^k\big\|_{L^1(D^k)}
\]
for every $t>0$.
In particular, there exists at most one entropy solution for given initial data.
\end{theorem}

\begin{proof}
From Lemma \ref{lem:locoinvariant} it follows that $\bu(t),\bv(t)\in\initialdataset(\germ)$ for a.e.~$t>0$. Let \( k\in\Indin \); the case \( k\in\Indout \) will follow analogously. The first step is a standard doubling of variables argument on each edge $k\in\Indin$ by selecting $\phi^k\in C_c^\infty(D^k\times[0,\infty))$ and $\phi_l\equiv0$ for $l\neq k$. The doubling of variables technique on a single edge gives:
\begin{equation}\label{eqn:doubledVariables}
\begin{split}
\int_{D^k} \int_0^{\infty} \big|u^k(x,t) - v^k(x,t)\big| \phi_t + q_v^k(u) \phi_x \,dt \,dx& \\
 + \int_{D^k} \big|\bar{u}^k(x) - \bar{v}^k(x)\big|\phi(x, 0) \,dx& \geq 0.
\end{split}
\end{equation}
Next, for general $\phi^k\in C_c^\infty\big(\overline{D^k}\times[0,\infty)\big)$, we cut off the functions near $x=0$ and couple the terms \eqref{eqn:doubledVariables} on each edge together by utilizing \eqref{eq:RHentr}. For $h>0$ we define
\[\mu_h (x) \coloneqq
\begin{cases}
0                            & x \in (-\infty,-2h) \\
\frac{1}{h}(x+2h) & x \in [-2h, -h) \\
1                            & x \in [-h, 0]
\end{cases}\]
and
\[ \Psi_h (x) \coloneqq 1 - \mu_h(x). \]
The derivative of \(\Psi_h \) reads
\[\Psi_h'(x) =
\begin{cases}
0           & x \in (-\infty,-2h) \\
-\frac{1}{h} & x \in [-2h, -h) \\
0           & x \in [-h, 0].
\end{cases}\]
Define \( \phi^k(x,t) \coloneqq \xi^k (x, t) \Psi_h (x) \) for a function \( \xi^k \in C_c^\infty \big( \overline{D^k} \times [0, \infty) \big) \).
We insert \( \phi \) into equation \eqref{eqn:doubledVariables} to get
\begin{align*}
\int_{D^k} \int_0^\infty \big|u^k(x,t) - v^k(x,t)\big| \xi_t^k \Psi_h + q^k_{v^k}(u^k) \xi_x^k \Psi_h \,dt \,dx& \\
+ \int_{D^k} \int_0^\infty q^k_{v^k}(u^k) \xi^k \Psi_h' \,dt\,dx  + \int_{D^k} \big|\bar{u}^k(x)-\bar{v}^k(x)\big| \xi^k \Psi_h \,dx& \geq 0.
\end{align*}
Sending \( h \downarrow 0 \) we get
\begin{align*}
\int_{D^k} \int_0^\infty  \big|u^k(x,t) - v^k(x,t)\big| \xi_t^k + q^k_{v^k}(u^k) \xi_x^k \,dt\,dx + \int_{D^k} \big|\bar{u}^k(x) - \bar{v}^k(x)\big| \xi^k \,dx & \\
+ \lim_{h \downarrow 0} \int_{-2h}^{-h} \int_0^\infty q^k_{v^k}(u^k) \xi^k \Psi_h' \,dt\,dx &\geq 0.
\end{align*}
Since the traces of \( q^k (u^k) \) and \( q^k (v^k) \) exist, we get
\begin{align*}
-\lim_{h \downarrow 0} \frac{1}{h} \int_0^T \int_{-2h}^{-h} q^k_{v^k}(u^k) \xi^k \,dx \,dt = -\int_0^T q_{v^{k-}}^k (u^k) \xi^k (0,t) \,dt.
\end{align*}
We therefore obtain
\begin{equation}\label{boundaryValue}
\begin{split}
\int_{D^k} \int_0^\infty \big|u^k(x,t) - v^k(x,t)\big| \xi_t + q^k_{v^k}(u^k) \xi_x \,dt\,dx& \\
+ \int_{D^k} \big|\bar{u}^k(x) - \bar{v}^k(x)\big| \,dx
-\int_0^T q_{v^k}^k (u^k) \xi(0,t) \,dt &\geq 0.
\end{split}
\end{equation}
By an analogous argument we get
\begin{align*}
\int_{D^k} \int_0^\infty \big|u^k(x,t) - v^k(x,t)\big| \xi_t + q^k_{v^k}(u^k) \xi_x \,dt\,dx& \\
+ \int_{D^k} \big|\bar{u}^k(x) - \bar{v}^k(x)\big| \,dx + \int_0^T q_{v^k}^k (u^k) \xi(0,t) \,dt& \geq 0
\end{align*}
for \( k\in\Indout\).
Fix $s>0$, let $r,\kappa>0$, and let \(\alpha_r:\R_-\to\R\) and \(\beta_\kappa\colon\R_+\to\R\) be given by
\[
\alpha_r (x) =
\begin{cases}
0 & x\in(-\infty,-r-1] \\
x+r+1 & x\in(-r-1,-r) \\
1 & x \in[-r,0)
\end{cases} \]
\[
\beta_\kappa(t)=
\begin{cases}
1 & t \in [0, s] \\
\frac{1}{\kappa} (\kappa + s - t) & t \in ( s, s+\kappa) \\
0 & t \in [s+\kappa, \infty).
\end{cases}
\]
Via a standard regularization argument one can check that \( \phi ( x, t) = \alpha_r (x) \beta_\kappa (t) \) is an admissible test function.
We compute the partial derivatives of \(\phi \):
\[
\phi_t(x,t) =
\begin{cases}
0                                          & t \in [ 0, s] \\
-\frac{1}{\kappa}\alpha_r (x) & t \in (s, s+\kappa] \\
0                                          & t \in (s + \kappa, \infty)
\end{cases}\]
and
\[
\phi_x(x,t) =
\begin{cases}
0                                                        & x\in(-\infty,-r-1) \\
\beta_\kappa (t) & x\in(-r-1,-r) \\
0                                                        & x\in(-r,0).
\end{cases}\]
We insert this into \eqref{boundaryValue} to get
\begin{align*}
 -\frac{1}{\kappa} \int_s^{s+\kappa}  \int_{-r-1}^0 \big|u^k(x,t) - v^k(x,t)\big| \alpha_r(x) \,dx \,dt& \\
 + \int_0^{s+\kappa} \int_{-r-1}^{-r} q_{v^k}^k (u^k) \beta_\kappa(t) \,dx\,dt + \int_{-r-1}^0 \big|\bar{u}^k(x) - \bar{v}^k(x)\big| \alpha_r (x) \,dx& \\
 - \int_0^{s+\kappa} q_{v^k}^k ( u^k (0, t)) \beta_\kappa (t) \,dt &\geq 0.
\end{align*}
Letting \( \kappa \to 0 \) and \( r\to\infty \), we get
\begin{equation*}
\big\| u^k(x,t) - v^k(x,t) \big\|_{L^1 ( D^k)} \leq \big\| \bar{u}^k(x) - \bar{v}^k(x)\big\|_{L^1 ( D^k)}  - \int_0^s q_{v^k}^k ( u^k(0,t)) \,dt.
\end{equation*}
An analogous inequality holds for $k\in\Indout$. We sum over all edges to get
\begin{align*}
&\sum_{k\in\Ind}  \big\|u^k(x,t) - v^k(x,t)\big\|_{L^1(D^k)} \\
&\leq \sum_{k\in\Ind} \int_{D^k} \big|\bar{u}^k(x) - \bar{v}^k(x)\big| \,dx +\underbrace{\int_0^s \sum_{k\in\Ind} \sgn(k) q^k_{v^k}(u^k)}_{\substack{\text{$\leq 0$ by \eqref{eq:RHentr} and}\\
\text{ maximality of $\germ$}}} \\
& \leq \sum_{k\in\Ind} \big\|\bar{u}^k(x) - \bar{v}^k(x)\big\|_{L^1(D^k)}.
\end{align*}
\end{proof}
%


%
\section{Numerical approximation}\label{Sec:Scheme}
In this section we construct a finite volume numerical approximation for \eqref{eq:clnetw} and prove stability and convergence properties of the method. The numerical method is rather standard for hyperbolic conservation laws, but an important feature of the method is that the vertex is discretized as a separate control volume. Although this control volume vanishes as the mesh parameter $\Dx$ is passed to zero, its presence will ensure that entropy is correctly dissipated at the vertex, even in the limit $\Dx\to0$.

\subsection{A finite volume method on networks}
Let $\Dt,\Dx>0$ be given discretization parameters. We define the index sets
\begin{align*}
D_\disc^+ \coloneqq \N, \qquad D_\disc^- \coloneqq -\N, \qquad D_\disc^k \coloneqq D_\disc^{\sgn(k)}, \qquad D_\disc^0 \coloneqq \{0\}.
\end{align*}
For $n\in\N_0$ we discretize\footnote{In numerical experiments, the timestep $\Dt$ is chosen dynamically at each step $n$ in order to comply with the CFL condition derived in Section~\ref{Sec:Scheme}. We use a uniform timestep for the sake of simplicity only.} $t^n=n\Dt$, and for $k\in\Ind$ and $i\in\Z$ we let $x_\iphf = (\iphf)\Dx$, and partition the physical domain into \emph{cells}
\[
\cell_i^k = D^k \cap \big(x_\imhf,\, x_\iphf\big).
\]
We define the mesh size at the vertex by $\Dx_0\coloneqq\sum_{k\in\Ind}|\cell_0^k|=N\Dx/2$, where $|A|$ denotes the Lebesgue measure of $A\subset\R$.
We make the finite volume approximation
\begin{align*}
u_i^{k,n} &\approx \frac{1}{\Dx}\int_{\cell_i^k}u^k(x,t^n)\,dx \qquad \text{for } i\in D^k_\disc, \\
u_0^n &\approx \frac{1}{\Dx_0}\sum_{k\in\Ind}\int_{\cell_0^k} u^k(x,t^n)\,dx.
\end{align*}
Fix some \( i \in D_\disc^k \), let $\phi^k(x,t) = \frac{1}{\Dt\Dx}\ind_{\cell_i^k}(x)\ind_{[t^n,\,t^{n+1})}(t)$ and $\phi^l\equiv0$ for $l\neq k$, and (after an approximation procedure) insert these into \eqref{eq:weakform}.
We then obtain the numerical method
\begin{subequations}\label{eq:fvm}
\begin{equation}
\frac{u_i^{k,n+1}-u_i^{k,n}}{\Dt} + \frac{F_\iphf^{k,n}-F_\imhf^{k,n}}{\Dx}=0
\end{equation}
where $F_\iphf^{k,n}=F^k\big(u^{k,n}_i,u^{k,n}_{i+1}\big)$ is an approximation of the mean flux through $x_\iphf$ over the time interval $[t^n,t^{n+1})$,
\[
F_\iphf^{k,n}\approx\frac{1}{\Dt}\int_{t^n}^{t^{n+1}}f^k\big(u^k(x_\iphf,t)\big)\,dt.
\]
For the special cell $i=0$ we let $\phi^k(x,t) = \frac{1}{\Dt\Dx_0}\ind_{\cell_0^k}(x)\ind_{[t^n,\,t^{n+1})}(t)$ for $k\in\Ind$ to obtain
\begin{equation}
\frac{u_0^{n+1}-u_0^n}{\Dt} + \frac{1}{\Dx_0}\biggl(\sum_{k\in\Indout} F_\hf^{k,n} - \sum_{k\in\Indin} F_{-\hf}^{k,n} \biggr) = 0.
\end{equation}
\end{subequations}
We will use the notational convention that $u_0^{k,n} \equiv u_0^n$ for all \(k\in\Ind\). (We postpone the definition of the initial data $u_i^{k,0}$ until Section~\ref{Sec:stabilityOfFVM}.)

Given a numerically computed solution $(\bu_i^n)_{i,n}$, we define the piecewise constant function
\begin{equation}\label{eq:proj}
\bu_\Dt(x,k,t) = u_i^{k,n} \qquad \text{for } x\in \cell_i^k,\ t\in[t^n,t^{n+1}).
\end{equation}
We remark that the integral of $\bu_\Dt$ w.r.t.~the measure $\lambda$ (cf.~\eqref{eq:measuredef}) can be written
\begin{equation}\label{eq:integralOfNumSoln}
\int_\Omega \bu_\Dt(\cdot,t)\,d\lambda = \sum_{k\in\Ind}\sum_{i\in D_\disc^k} u_i^{k,n}\,\Dx + u_0^n\,\Dx_0
\end{equation}
for any $t\in[t^n,t^{n+1})$, and the total variation of $\bu_\Dt$ (cf.~\eqref{eq:tvdef}) can be written
\begin{equation}\label{eq:tvdiscrete}
\begin{split}
\TV(\bu_\Dt(\cdot,t)) =& \sum_{k\in\Indin}\sum_{i \in D_\disc^k}\big|u_{i+1}^{k,n}-u_i^{k,n}\big| + \sum_{k\in\Indout}\sum_{i \in D_\disc^k}\big|u_{i}^{k,n}-u_{i-1}^{k,n}\big| \\
=& \sum_{k\in\Indin}\sum_{i \in D_\disc^k}\big|u_{i}^{k,n}-u_{i-1}^{k,n}\big| + \sum_{k\in\Indout}\sum_{i \in D_\disc^k}\big|u_{i+1}^{k,n}-u_{i}^{k,n}\big| \\
&+ \sum_{k\in\Indin}\big|u_0^{n}-u_{-1}^{k,n}\big| + \sum_{k\in\Indout}\big|u_0^{n}-u_{1}^{k,n}\big|.
\end{split}
\end{equation}
Note also that a numerical method of the form \eqref{eq:fvm} is
\emph{conservative} in the sense that the total mass $\sum_{k\in\Ind}
\int_\Omega \bu_\Dt \,d\lambda$ is independent of $n$:
\begin{align*}
\int_\Omega \bu_\Dt(\cdot,t^{n+1}) \,d\lambda ={}& \sum_{k\in\Ind} \sum_{i \in D_\disc^k} u_i^{k,n+1} \Dx + u_0^{n+1} \Dx_0 \\
={}&\sum_{k\in\Ind} \sum_{i \in D_\disc^k} u_i^{k,n} \Dx - \Dt \big(F_\iphf^{k,n} - F_\imhf^{k,n}\big) + u_0^n \Dx_0 \\
&{{}-\Dt} \biggl( \sum_{k\in\Indout} F_{\hf}^k  - \sum_{k\in\Indin} F_{-\hf}^k \biggr) \\
={}& \sum_{k\in\Ind} \sum_{i\in D_\disc^k} u_i^{k,n} \Dx + u_0^n \Dx_0 \\
={}& \int_\Omega \bu_\Dt(\cdot,t^{n}) \,d\lambda.
\end{align*}

As a shorthand for the scheme \eqref{eq:fvm} we define the functions
\begin{subequations}\label{eq:UpdateFct}
\begin{equation}\label{eq:UpdateFctEdge}
G^k \big( u_{i-1}^k, u_i^k, u_{i+1}^k \bigr) \coloneqq u_i^{k} - \frac{\Dt}{\Dx}\big(F^k\bigl(u_i^k, u_{i+1}^k\bigr) - F^k \bigl( u_{i-1}^k, u_i^k\bigr)\big)
\end{equation}
for $k\in\Ind$ and
\begin{equation}\label{eq:UpdateFctNode}
\begin{split}
  &G^0\bigl( u_{-1}^{-\Nin}, \dots, u_{-1}^{-1}, u_0, u_1^1, \dots,
  u_1^{\Nout}\bigr) \\
  &\quad\coloneqq
  u_0 - \frac{\Dt}{\Dx_0} \biggl(\sum_{k\in\Indout} F^k\bigl( u_0, u_1^k\bigr)
  - \sum_{k\in\Indin} F^k\bigl(u_{-1}^k, u_0\bigr)\biggr),
\end{split}
\end{equation}
\end{subequations}
enabling us to write \eqref{eq:fvm} in the update form
\begin{equation}\label{eq:updateform}
\begin{split}
u_i^{k,n+1} &= G^k \big( u_{i-1}^{k,n}, u_i^{k,n}, u_{i+1}^{k,n} \bigr) \qquad\qquad \text{for } i\in D_\disc^k,\ k\in\Ind \\
u_0^{n+1} &= G^0\bigl( u_{-1}^{-\Nin,n}, \dots, u_{-1}^{-1,n}, u_0^n, u_1^{1,n}, \dots, u_1^{\Nout,n}\bigr).
\end{split}
\end{equation}
As a shorthand for \eqref{eq:updateform}, we will sometimes use the notation
\begin{equation}\label{eq:updateformvec}\tag{\ref{eq:updateform}'}
u_i^{k,n+1} = G^k\big(\bu_{i-1}^n, \bu_i^n, \bu_{i+1}^n\big) \qquad \text{for }i\in D_\disc^k,\ k\in\Ind_0,
\end{equation}
where $\bu_i^n$ is the vector containing all numerical values at index $i$ at time $n$.



\begin{definition}[Monotone scheme]
The difference scheme \eqref{eq:updateformvec} is \emph{monotone} if
\[
\bu^n\leq\bv^n \quad \Rightarrow \quad \bu^{n+1}\leq \bv^{n+1},
\]
where $\bu^n\leq\bv^n$ means that every component $u_i^{k,n}$ of $\bu^n$ is no greater than the corresponding component of $\bv^n$.
\end{definition}
We state a straightforward CFL-type condition which ensures monotonicity of the numerical scheme.
\begin{proposition}
Consider a consistent finite volume method \eqref{eq:fvm}, where $F^k$ is nondecreasing in the first variable and nonincreasing in the second. Then the method is monotone under the CFL condition
\begin{equation}\label{eq:cfl}
\Dt\max_{k,u,v}\Bigl|\frac{\partial F^k}{\partial u}(u,v)\Bigr| \leq \Dx/2, \qquad \Dt\max_{k,u,v}\Bigl|\frac{\partial F^k}{\partial v}(u,v)\Bigr| \leq \Dx/2.
\end{equation}
\end{proposition}

\begin{proof}
We can calculate the derivatives to the update functions to get
\begin{align*}
  \frac{\partial G^k}{\partial u_{i-1}^{k}}
  &= \frac{\Dt}{\Dx}
    \frac{\partial
    F_\imhf^k}{\partial
    u_{i-1}^{k}}, \qquad \frac{\partial G^k}{\partial u_{i+1}^{k}} =
    -\frac{\Dt}{\Dx} \frac{\partial F_\iphf^k}{\partial u_{i+1}^{k}}
  \\
  \frac{\partial G^k}{\partial u_{i}^{k}}
  &= 1 - \frac{\Dt}{\Dx} \biggl( \frac{\partial F_\iphf^k}{\partial u_{i}^{k}} -
    \frac{\partial F_\imhf^k}{\partial u_{i}^{k}} \biggr),
\end{align*}
for each $k\in\Ind$, and
\begin{align*}
  \frac{\partial G^0}{\partial u_{-1}^{k}}
  &= \frac{\Dt}{\Dx_0} \frac{\partial F_{-\hf}^k}{\partial u_{-1}^{k}} && \text{for } k\in\Indin, \\
  \frac{\partial G^0}{\partial u_{1}^{k}}
  &= -\frac{\Dt}{\Dx_0} \frac{\partial F_{\hf}^k}{\partial u_{1}^{k}} && \text{for } k\in\Indout, \\
  \frac{\partial G^0}{\partial u_0}
  &= 1 - \frac{\Dt}{\Dx_0} \biggl(\sum_{k\in\Indout}\frac{\partial F_{\hf}^k}{\partial u_0^n} - \sum_{k\in\Indin}\frac{\partial F_{-\hf}^k}{\partial u_0^n} \biggr)
\end{align*}
on the vertex. We would like these derivatives to be non-negative. The monotonicity of $F^k$ guarantees that the first, second, fourth and fifth expressions are non-negative. Applying monotonicity of $F^k$ to the third and sixth terms, we get
\[
\frac{\partial G^k}{\partial u_{i}^{k}} = \frac{\Dt}{\Dx}
\biggl(\frac{\Dx}{\Dt} - \Bigl|\frac{\partial F_\iphf^k}{\partial
      u_{i}^{k}}\Bigr| - \Bigl|\frac{\partial F_\imhf^k}{\partial u_{i}^{k}} \Bigr|\biggr) \geq 0
\]
(by \eqref{eq:cfl}) and
\begin{align*}
\frac{\partial G^0}{\partial u_0} &= \frac{\Dt}{\Dx_0}\biggl(\frac{\Dx_0}{\Dt} - \sum_{k\in\Indout}\Bigl|\frac{\partial F_{\hf}^k}{\partial u_0^n}\Bigr| - \sum_{k\in\Indin}\Bigl|\frac{\partial F_{-\hf}^k}{\partial u_0^n}\Bigr|\biggr) \\
\intertext{(using $\Dx_0=N\Dx/2$)}
&\geq \frac{\Dt}{\Dx_0}\biggl(N\frac{\Dx}{2\Dt} - \Nout\max_{k,u,v}\Bigl|\frac{\partial F^k}{\partial u}(u,v)\Bigr| - \Nin\max_{k,u,v}\Bigl|\frac{\partial F^k}{\partial v}(u,v)\Bigr|\biggr)\\
&\geq 0
\end{align*}
by \eqref{eq:cfl}.
\end{proof}

\subsection{Discrete stationary solutions}
In the same way that stationary solutions are essential for the well-posedness of entropy solutions (cf.~Section~\ref{Sec:Stability}), they are essential to the stability and convergence of numerical methods on networks. Asserting that a numerical solution is both constant in time and on each edge yields the following definition.
\begin{definition}[Discrete Stationary Solution]
Consider a consistent, conservative numerical method \eqref{eq:fvm}. A \textit{discrete stationary solution} for \eqref{eq:fvm} is a vector
\[ \bc^\disc \coloneqq
( c^{-\Nin}, \dots, c^{\Nout}) \in\R^{N+1}
\]
satisfying the Rankine--Hugoniot condition
\begin{equation}\label{eq:RHdss}
\sum_{k\in\Indin} f^k ( c^k ) = \sum_{k\in\Indout} f^k ( c^k )
\end{equation}
as well as the conditions
\begin{subequations}\label{eq:dssEdge}
\begin{align}
F^k ( c^k, c^0 ) &= f^k ( c^k ) \qquad \text{for } k\in\Indin,  \label{eq:consistencyInEdge} \\
F^k ( c^0, c^k ) &= f^k ( c^k ) \qquad \text{for }k\in\Indout. \label{eq:consistencyOutEdge}
\end{align}
\end{subequations}
\end{definition}

In the remainder, sets of discrete stationary solutions will be denoted with a superscript, $\germ^0$, to signal that they also include a value at the vertex $i=0$.

\begin{notation}\label{not:dssindex}
We will sometimes index a discrete stationary solution as
\begin{subequations}\label{eq:dssindex}
\begin{equation}\label{eq:dssindex1}
\bc_i = \begin{cases}
\big(c^{-\Nin},\dots,c^{-1}\big) & i<0 \\
c^0 & i=0 \\
\big(c^1,\dots,c^\Nout\big) & i>0
\end{cases}
\end{equation}
for $i\in\Z$ and, by extension,
\begin{equation}\label{eq:dssindex2}
c_i^k = \begin{cases}
c^k & i\neq 0 \\
c^0 & i=0.
\end{cases}
\end{equation}
\end{subequations}
Using the notation \eqref{eq:updateformvec}, it is readily checked that discrete stationary solutions are precisely those that are constant on each edge and satisfy
\[
\bc_i = G^k(\bc_{i-1},\bc_i,\bc_{i+1}) \qquad \forall\ i\in D_\disc^k,\ k\in\Ind_0.
\]
\end{notation}

\begin{remark}
The conditions \eqref{eq:dssEdge} say that the numerical fluxes at the vertex reduce to the upwind flux on the in edges and the downwind flux on the out edges. This can be interpreted as information only flowing into the vertex, not out of it. This is consistent with the interpretation of the vertex as a stationary shock.
\end{remark}

\begin{remark}
Discrete stationary solutions \( \bc = \big(c^{-\Nin}, \dots, c^{\Nout}\big) \) fulfil a discrete version of the Rankine--Hugoniot type condition \eqref{eqn:RankingeHugoniot},
\[ \sum_{k\in\Indin} F^k\big( c^k, c^0\big) = \sum_{k\in\Indout} F^k\big( c^0, c^k\big). \]
\end{remark}

\begin{lemma}\label{lem:CharacterizationDiscreteStationarySol}
Consider a consistent, conservative numerical scheme \eqref{eq:fvm}. Let \( \bc = ( c^k )_{k\in\Ind} \) be a stationary solution for \eqref{eq:cl} and let $c^0\in\R$. Then the vector \( \bc_\disc = ( c^{-\Nin}, \dots, c^{-1}, c^0, c^1, \dots, c^{\Nout} )\in\R^{N+1} \) is a discrete stationary solution if and only if
\[
c^0 \in  \bigcap_{k\in\Indin} ( H^k )^{-1} (\{f^k(c^k)\}) \mathrel{\bigcap}\bigcap_{k\in\Indout} ( J^k )^{-1} (\{f^k(c^k)\})
\]
where
\begin{align*}
H^k(c) &\coloneqq F^k(c^k,c) \quad \text{for } k\in\Indin, \qquad J^k(c) \coloneqq F^k(c,c^k) \quad \text{for } k\in\Indout.
\end{align*}
\end{lemma}
\begin{proof}
We can rewrite conditions \eqref{eq:consistencyInEdge} and \eqref{eq:consistencyOutEdge} as
\begin{align*}
&\eqref{eq:consistencyOutEdge} \quad \Leftrightarrow \quad H^k(c^0)=f^k(c^k) \quad \Leftrightarrow \quad c^0 \in ( H^k )^{-1} (\{f^k (c^k)\})
\intertext{for $k\in\Indin$, and} \\
&\eqref{eq:consistencyInEdge} \quad \Leftrightarrow \quad J^k(c^0)=f^k(c^k) \quad \Leftrightarrow \quad c^0 \in ( J^k )^{-1} (\{f^k (c^k)\})
\end{align*}
for $k\in\Indout$. Hence, if \eqref{eq:consistencyInEdge}, \eqref{eq:consistencyOutEdge} are satisfied then $c^0$ must lie in all of the sets on the right hand side, and hence in their intersection. Conversely, if $c^0$ lies in the intersection, then \eqref{eq:consistencyInEdge}, \eqref{eq:consistencyOutEdge} are satisfied.
\end{proof}

\subsection{A stability framework for general consistent, conservative, monotone methods}\label{Sec:stabilityOfFVM}
We set out to prove an \emph{\(L^\infty\) bound}, \textit{\(L^1\) contractiveness} and \emph{Lipschitz continuity in time} for solutions computed with a general consistent, conservative, monotone finite volume method on a network. Our starting point will be a class of discrete stationary solutions $\germ_\disc^0\subset\R^{N+1}$ for a conservative finite volume method \eqref{eq:fvm}. We take initial data $\bar u\in \initialdataset(\germ_\disc^0)$ (cf.~\eqref{eq:definitialdataset}), we let $\bc\in\germ_\disc^0$ be as specified in \eqref{eq:definitialdataset}, and consider the finite volume method \eqref{eq:fvm} initialized by
\begin{equation}\label{eq:fvminitialdata}
u_i^{k,0} = \frac{1}{\Dx}\int_{\cell_i^k}\bar u^k(x)\,dx, \qquad u_0^0 = c^0.
\end{equation}
(The value $c^0$ is chosen for convenience, and any value in $[c^0,d^0]$ will have the desired effect.)

\begin{lemma}
Consider monotone numerical flux functions $F^k$ ($k\in\Ind$). Let \(\bc, \bd\) be discrete stationary solutions satisfying $c^k\leq d^l$ for all $k,l\in\Ind$ (cf.~Definition~\ref{Def:Loco}). Then $c^0,d^0$ can be modified such that $\bc,\bd$ remain discrete stationary solutions and such that $c^0\leq d^0$.
\end{lemma}
\begin{proof}
Define
\[
I^k(c^k) \coloneqq \begin{cases} \big(F^k(c^k,\cdot)\big)^{-1} \big(\{ f^k(c^k) \}\big) & \text{for } k \in \Indin \\
\big(F^k(\cdot,c^k)\big)^{-1}\big(\{ f^k (c^k) \}\big) &  \text{for } k \in \Indout.
\end{cases}
\]
Since all \( F^k \) are monotone, each \( I^k(c^k) \) is a connected interval which contains \( c^k \), and moreover, Lemma~\ref{lem:CharacterizationDiscreteStationarySol} says that \( c^0 \in \bigcap_{k\in\Ind} I^k(c^k) \). This implies that $\convhull{c^0}{c^k} \subset I^k(c^k)$, where $\convhull{a}{b}=[\min(a,b),\,\max(a,b)]$. Since $\bigcap_{k\in\Ind}\convhull{c^0}{c^k}$ is nonempty, the number
\[
\tilde{c}^0 \coloneqq \min \biggl(\bigcap_{k\in\Ind}\convhull{c^0}{c^k}\biggr)
\]
exists and satisfies $\tilde{c}^0\leq \max_{k\in\Ind} c^k$. Appealing again to Lemma~\ref{lem:CharacterizationDiscreteStationarySol}, $\bc$ remains a discrete stationary solution if $c^0$ is replaced by $\tilde{c}^0$. In a similar way we replace $d^0$ by
\[
\tilde{d}^0 \coloneqq \max \biggl(\bigcap_{k\in\Ind}\convhull{d^0}{d^k}\biggr),
\]
which satisfies $\tilde{d}^0\geq \min_{k\in\Ind} d^k$. By our hypothesis, it follows that $\tilde c^0 \leq \tilde{d}^0$.
\end{proof}
\begin{proposition}\label{prop:linfbound}
Consider a consistent, conservative, monotone finite volume method \eqref{eq:fvm}, \eqref{eq:fvminitialdata} with a set of discrete stationary states $\germ_\disc^0$. For any initial data $\bar \bu \in \initialdataset(\germ_\disc^0)$, the numerical solution is uniformly $L^\infty$ bounded.
\end{proposition}
\begin{proof}
Pick discrete stationary states $\bc,\bd\in\germ_\disc^0$ as in \eqref{eq:definitialdataset}. It is clear that the initial data defined in \eqref{eq:fvminitialdata} satisfy $c^k \leq u_i^{k,0} \leq d^k$ for all $i\in D_\disc^k$ and $k\in\Ind_0$. If the same holds at some time step $n\in\N_0$ then (using the notation \eqref{eq:updateformvec}, \eqref{eq:dssindex})
\[
u_i^{k,n+1} = G^k \bigl( \bu_{i-1}^{n}, \bu_i^{n}, \bu_{i+1}^{n} \bigr) \geq G^k\bigl( \bc_{i-1}, \bc_i, \bc_{i+1}\bigr) = c_i^k
\]
for all $i\in D_\disc^k$ and $k\in\Ind_0$, and similarly, $u_i^{k,n+1} \leq d_i^k$.
\end{proof}
\begin{definition}[\(L^1\) contractive method]
A numerical method \eqref{eq:updateformvec} is \emph{\( L^1 \) contractive} if
\begin{align*}
\bigl\| \bu_\Dt ( \cdot, t ) - \bv_\Dt ( \cdot, t ) \bigr\|_{L^1(\Omega;\lambda)} \leq \|\bar \bu - \bar \bv\|_{L^1(\Omega;\lambda)}
\end{align*}
for all \( t \geq 0 \), where $\bu_\Dt$ and $\bv_\Dt$ are the projection of the numerical solution (cf.~\eqref{eq:proj}) computed with initial data $\bar \bu,\bar \bv\in \initialdataset(\germ_\disc^0)\cap L^1(\Omega;\lambda)$, respectively. (See~\eqref{eq:integralOfNumSoln} for the integral of $\bu_\Dt$, $\bv_\Dt$ w.r.t.~$\lambda$.)
\end{definition}
We state the well known Crandall--Tartar lemma which we will use in the following proof. Here and below, we use the notation $a\vee b = \max(a,b)$.
%
%
\begin{theorem}[Crandall--Tartar: {\cite[Proposition 1]{CrandallTartar1980}}]\label{CrandallTartar}
Let \( (\Omega,\lambda) \) be a measure space. Let \( C \subset L^1 ( \Omega;\lambda ) \) have the property that \( f,g \in C \) implies \( f \vee g \in C \). Let \( V\colon C \to L^1(\Omega;\lambda) \) satisfy \( \int_\Omega V ( f )\,d\lambda = \int_\Omega f\,d\lambda \) for \( f \in C \).
Then the following three properties of \( V \) are equivalent:
\begin{itemize}
\item[(a)] \( f, g \in C \) and \( f \leq g \) a.e. implies \( V( f) \leq V( g) \) a.e.,
\item[(b)] \( \int_\Omega ( V( f) - V( g) )^+ \leq \int_\Omega ( f - g )^+ \) for \( f, g \in C \),
\item[(c)] \( \int_\Omega \big|V( f) - V( g)\big| \leq \int_\Omega |f - g| \) for \( f, g \in C \).
\end{itemize}
\end{theorem}
%
%
%
We can now prove $L^1$-contractivity of our method.
\begin{theorem}
Every conservative, consistent monotone method \eqref{eq:fvm}, \eqref{eq:fvminitialdata} is \(L^1\)-contractive.
\end{theorem}
\begin{proof}
Let $C=\constfuncs$ be the set of piecewise constant functions,
\[
\constfuncs = \Bigl\{\bu\in L^1\cap L^\infty(\Omega;\lambda)\ :\ \bu(x)=\sum_{k\in\Ind_0}\sum_{i\in D_\disc^k}u_i^k \ind_{\cell_i^k} \text{ for } u_i^k\in\R \Bigr\}.
\]
We define the operator \( V:\constfuncs\to \constfuncs \) mapping a numerical solution to the next time step,
\begin{align*}
  V ( \bu )
  \coloneqq{}& \sum_{k\in\Ind} \sum_{i\in D_\disc^k} \ind_{\cell_i^k} \Bigl( u_i^{k} - \frac{\Dt}{\Dx} \big(F^k\big(u_i^k,u_{i+1}^k\big) - F^k\big(u_{i-1}^k,u_i^k\big)\big) \Bigr) \\
  &+ \sum_{k\in\Ind}\ind_{\cell_0^k} \biggl(u_0^{0} - \frac{\Dt}{\Dx_0} \biggl(
    \sum_{k\in\Indout} F^k(u_0,u_1^k) - \sum_{k\in\Indin}
    F^k(u_{-1}^k,u_0)
    \biggr) \biggr).
\end{align*}
By the definition \eqref{eq:measuredef} of the measure $\lambda$ (cf.\ also \eqref{eq:integralOfNumSoln}), we have \( \int_{\Omega} V ( \bu ) \,d\lambda = \int_{\Omega} \bu \,d\lambda \) for all $\bu\in \constfuncs$. We apply the Crandall--Tartar lemma to conclude $L^1$-contractivity of the numerical solution operator $V$.
\end{proof}

From \( L^1 \)-contractivity we get continuity in time as a corollary:

\begin{corollary}\label{cor:timeContinuity}
Consider a consistent, conservative and monotone method
\eqref{eq:fvm}. Let $\bu_\Dt$ be an approximate solution computed with
this method and assume that all numerical fluxes $F^k$ are Lipschitz
continuous in both arguments. Then computed solutions are uniformly
$L^1$ Lipschitz continuous in time, i.e.,
\begin{align*}
\big\|\bu_\Dt(t^{n+1}) - \bu_\Dt(t^n) \big\|_{L^1(\Omega;\lambda)} & \leq \big\|\bu_\Dt(t^1) - \bu_\Dt(t^0) \big\|_{L^1(\Omega;\lambda)} \\
& \leq \Dt \big(C\TV(\bu^0) + \bar M\big),
\end{align*}
where the constants $C$ and $\bar{M}$ depend on the flux functions and
on the initial data.
\end{corollary}
\begin{proof}
We compute
\begin{align*}
  \bigl\|&\bu_\Dt(t^{n+1}) - \bu_\Dt(t^n) \bigr\|_{L^1(\Omega;\lambda)} \\
&= \big\|V(\bu_\Dt(t^{n})) - V(\bu_\Dt(t^{n-1})) \big\|_{L^1(\Omega;\lambda)} \\
\intertext{({using Theorem~\ref{CrandallTartar}(c)})}
&\leq \big\|\bu_\Dt(t^{n}) - \bu_\Dt(t^{n-1}) \big\|_{L^1(\Omega;\lambda)} \leq \dots \leq \big\|\bu_\Dt(t^1) - \bu_\Dt(t^0) \big\|_{L^1(\Omega;\lambda)} \\
&= \Dx \sum_{k\in\Ind} \sum_{i \in D_\disc^k} | u_i^{k, 1} - u_i^0| + \Dx_0|u_0^{k, 1} - u_0^0|  \\
&= \Dt \sum_{k\in\Ind} \sum_{i \in D_\disc^k} \big|F_\iphf^{k,0} - F_\imhf^{k,0}\big| + \Dt \biggl|\sum_{k\in\Indout} F_\hf^{k,0} - \sum_{k\in\Indin} F_{-\hf}^{k,0}\biggr| \\
&= \Dx \lambda \sum_{\in\Ind} \sum_{i\in D_\disc^k}\big|u_i^{k,0} - u_{i-1}^{k,0}\big| \\
&\quad+ \Dt \biggl|\!\begin{aligned}[t]&\sum_{k\in\Indout} F_\hf^{k,0} - F^{k,0} \big(u_0^0, u_0^0\big) - \sum_{k\in\Indin} F_{-\hf}^{k,0} - F^{k,0} (u_0^0, u_0^0) \\
&+\overbrace{\sum_{k\in\Indout} f^k \big(u_0^0\big)}^{\eqqcolon f_\out(u_0^0)} - \overbrace{\sum_{k\in\Indin} f^k\big(u_0^0\big)}^{\eqqcolon f_\into(u_0^0)}\biggr|\end{aligned} \\
&\leq \Dt \sum_{k\in\Ind} \sum_{i\in D_\disc^k} L^k \bigl(\big|u_i^{k,0} - u_{i-1}^{k,0}\big| + \big|u_{i+1}^{k,0} - u_i^{k,0}\big| \bigr) \\
&\quad + \Dt \underbrace{\Biggl(\sum_{k\in\Indout} L^k \big|u_1^{k,0} - u_0^0\big| + \sum_{k\in\Indin} L^k \big|u_0^0 - u_{-1}^{k,0}\big| + \overbrace{\big|f_\out\big(u_0^0\big) - f_\into\big(u_0^0\big)\big|}^{\leq M}\Biggr)}_{\eqqcolon\bar M} \\
&\leq \Dt (C \TV(\bu^0) + \bar M ),
\end{align*}
where we collect all constants into the global constant \( C \). We can bound \( \big|f_\out \big(u_0^0\big) - f_\into \big(u_0^0\big)\big| \leq \bar M \) with a constant \( \bar M \in \R \) since \( f_\into, f_\out \) are continuous and \( u_0^0 \in L^\infty \).
\end{proof}

\section{Convergence of finite volume schemes}\label{Sec:Convergence}
We are now in place to prove convergence in the case where the flux
functions $f^k$ are strictly monotone. We do this by using the upwind method
where the numerical flux functions are defined by
\[
F^k(u,v)=\begin{cases} f^k(u) & \text{if
$f^k$ is increasing,} \\
f^k(v) & \text{if $f^k$ is decreasing.}
\end{cases}
\]
We shall
show that the set of discrete approximations is compact in
$L^\infty\big([0,\infty);L^1(\Omega;\lambda)\big)$, and that any limit limit is an
entropy solution. In particular, this convergence result establishes
existence of an entropy solution.
We show convergence to a weak solution by proving a Lax--Wendroff type theorem:
\begin{theorem}[Lax--Wendroff theorem]\label{thm:laxwendroff}
  Fix \( T > 0 \). Assume that each flux function \( f^k \) is locally
  Lipschitz continuous and strictly monotone. Let $\germ_\disc^0$ be a
  class of discrete stationary solutions for the upwind method and let \( \bu_\Dt \) be computed from the upwind method with initial data
  $\bar \bu \in \initialdataset(\germ_\disc^0) \cap L^1(\Omega;\lambda)$.
  Consider a subsequence \( \big(\bu_{\Dt_{\ell}}\big)_{\ell\in\N} \) such that
  \( \Dt_{\ell} \to 0 \) and
  \( \bu_{\Dt_{\ell}}\to \bu\) in \(L^\infty([0,T];L^1(\Omega;\lambda)) \) as
  $k\to\infty$.  Then the limit $\bu$ is the unique entropy solution to
  \eqref{eq:clnetw}, that is, $\bu$ satisfies \eqref{eq:entropysoln}.
\end{theorem}
\begin{proof}
We write $\Dx$ and $\Dt$ rather than  $\Dx_{\ell}$ and $\Dt_{\ell}$, and
we shall show that $u$ satisfies the entropy condition \eqref{eq:entropysoln} for every $\bc\in\germ_\disc^0$. Choosing stationary solutions $\bc,\bd\in\germ_\disc^0$ such that $\bc\leq \bu_\Dt\leq \bd$ (cf.~Proposition~\ref{prop:linfbound}) in particular shows that $u$ is a weak solution.

Let $\bc\in\germ_\disc^0$ and consider the Crandall--Majda numerical entropy fluxes
\begin{align*}
Q_\iphf^{k,n} = F^k \bigl( u_i^{k,n} \vee c^k_i, u_{i+1}^{k,n} \vee c^k_i \bigr) - F^k \bigl( u_i^{k,n} \wedge c^k, u_{i+1}^{k,n} \wedge c^k \bigr)
\end{align*}
for $i=0,1,\dots$ when $k\in\Indout$, and for $i=\dots,-2,-1$ when $k\in\Indin$, and
\begin{align*}
Q_{-\hf}^n = \sum_{k\in\Indin} Q_{-\hf}^{k,n}, \qquad Q_{\hf}^n = \sum_{k\in\Indout} Q_{\hf}^{k,n}
\end{align*}
(cf.~Notation~\ref{not:dssindex} for the definition of $c_i^k$). Recalling the definition  \eqref{eq:UpdateFct} of the update functions $G^k$, we see that
\begin{align*}
&G^k \bigl( u_{i-1}^{k,n} \vee c^k, u_i^{k,n} \vee c^k , u_{i + 1}^{k,n} \vee c^k \bigr) -
G^k \bigl( u_{i-1}^{k,n} \wedge c^k, u_i^{k,n} \wedge c^k , u_{i + 1}^{k,n} \wedge c^k \bigr)\\
&\qquad = \bigl|u_i^{k,n} - c^k\bigr| - \frac{\Dt}{\Dx} \bigl( Q_\iphf^k - Q_\imhf^k \bigr),
\end{align*}
for $k\in\Ind$ and $i \in D_\disc^k$. Hence,
\begin{equation}\label{eq:DiscEntropyIneqEdge}
\begin{split}
\big|u_i^{k,n+1}-c^k\big| &= u_i^{k,n+1}\vee c^k - u_i^{k,n+1}\wedge c^k \\
&= G^k \bigl(u_{i-1}^{k,n}, u_i^{k,n}, u_{i + 1}^{k,n}\bigr)\vee c^k - G^k \bigl(u_{i-1}^{k,n}, u_i^{k,n}, u_{i + 1}^{k,n}\bigr)\wedge c^k \\
&\leq G^k \bigl( u_{i-1}^{k,n} \vee c^k, u_i^{k,n} \vee c^k , u_{i + 1}^{k,n} \vee c^k \bigr) \\
&\quad- G^k \bigl( u_{i-1}^{k,n} \wedge c^k, u_i^{k,n} \wedge c^k , u_{i + 1}^{k,n} \wedge c^k \bigr) \\
&= \bigl|u_i^{k,n} - c^k\bigr| - \frac{\Dt}{\Dx} \bigl( Q_\iphf^{k,n} - Q_\imhf^{k,n} \bigr).
\end{split}
\end{equation}
Similarly, we find that
\begin{equation}\label{eq:DiscEntropyIneqVertex}
\bigl|u_0^{n+1} - c^0\bigr| - \bigl|u_0^n - c^0\bigr| + \frac{\Dt}{\Dx_0} \bigl( Q_\hf^n - Q_{-\hf}^n \bigr) \leq 0
\end{equation}
We choose \( T = M \Dt \) for a natural number \( M \), multiply the above $N+1$ inequalities with a test function $\phi$ and sum up to get
\begin{align*}
0&\geq\sum_{n=0}^{M}\sum_{k\in\Ind} \sum_{i \in D_\disc^k} \phi_i^{k,n}\biggl( \bigl( \big|u_i^{k,n+1} - c^k\big| - \big|u_i^{k,n} - c^k\big| \bigr) + \frac{\Dt}{\Dx} \bigl( Q_\iphf^{k,n} - Q_\imhf^{k,n} \bigr)\biggr) \\
&\quad +\sum_{n=0}^{M} \phi_0^n\biggl( \frac{N}{2} \bigl( \big|u_0^{n+1} - c^0\big| - \big|u_0^n - c^0\big| \bigr) + \frac{\Dt}{\Dx} \biggl( \sum_{k\in\Indout} Q_{\hf}^{k,n} - \sum_{k\in\Indin} Q_{-\hf}^{k,n} \biggr) \biggr),
\end{align*}
where $\phi_i^{k,n}=\phi^k(x_i,t^n)$. After summation by parts we get
\begin{align*}
0&\geq-\sum_{n=1}^{M} \sum_{k\in\Ind} \sum_{i \in D_\disc^k} \big|u_i^{k,n} - c^k\big| \bigl( \phi_i^{k,n} - \phi_i^{k,n-1} \bigr) - \sum_{k\in\Ind} \sum_{i \in D_\disc^k} \big|u_i^{k,0} - c^k\big| \phi_i^{k,0} \\
&\quad -\frac{\Dt}{\Dx} \sum_{n=0}^{M} \biggl(\begin{aligned}[t]& \sum_{k\in\Indin} \sum_{i \in D_\disc^k} Q_\imhf^{k,n} \bigl( \phi_i^{k,n} - \phi_{i-1}^{k,n} \bigr) \\
&+ \sum_{k\in\Indout} \sum_{i \in D_\disc^k} Q_\iphf^{k,n} \bigl( \phi_{i+1}^{k,n} - \phi_i^{k,n} \bigr) \\
&+ \sum_{k\in\Indin} \phi_{-1}^{k,n} Q_{-\hf}^{k,n} - \sum_{k\in\Indout} \phi_1^{k,n} Q_{\hf}^{k,n} \biggr)\end{aligned} \\
&\quad-\frac{N}{2} \big|u_0^n - c^0\big| \phi_0^0 - \frac{N}{2}
\sum_{n=1}^{M} \big|u_0^n - c^0\big| \bigl( \phi_0^n - \phi_0^{n-1}
\bigr) \\
&\quad+ \frac{\Dt}{\Dx} \sum_{n=0}^{M} \biggl(\sum_{k\in\Indout} Q_{\hf}^n \phi_0^n - \sum_{k\in\Indin} Q_{-\hf}^n \phi_0^n \biggr).
\end{align*}
After shifting the $i$ index on the second line we get
\begin{align*}
0&\geq-\underbrace{\Dt \Dx \sum_{n=1}^{M} \sum_{k\in\Ind} \sum_{i \in D_\disc^k} \big|u_i^{k,n} - c^k\big| \biggl( \frac{\phi_i^{k,n} - \phi_i^{k,n-1}}{\Dt} \biggr)}_{=\,A_1} \\
&\quad- \underbrace{\Dx \sum_{k\in\Ind} \sum_{i \in D_\disc^k} \big|u_i^{k,0} - c^k\big| \phi_i^{k,0}}_{=\,A_2} \\
&\quad- \underbrace{\Dt \Dx \sum_{n=0}^{M} \biggl(\begin{aligned}[t]&\sum_{k\in\Indin} \sum_{i \in D_\disc^k} Q_\iphf^{k,n} \biggl( \frac{\phi_{i+1}^{k,n} - \phi_{i}^{k,n}}{\Dx} \biggr) \\
&+ \sum_{k\in\Indout} \sum_{i \in D_\disc^k} Q_\imhf^{k,n} \biggl( \frac{\phi_{i}^{k,n} - \phi_{i-1}^{k,n}}{\Dx} \biggr)\biggr)\end{aligned}}_{=\,A_3} \\
&\quad- \underbrace{\Dx \frac{N}{2} \big|u_0^{0} - c^0\big| \phi_0^0 - \Dt \Dx \frac{N}{2} \sum_{n=1}^{M} \big|u_0^n - c^0\big| \biggl( \frac{\phi_0^n - \phi_0^{n-1}}{\Dt} \biggr)}_{=\,A_4}.
\end{align*}
Taking limits we get for \( \Dt, \Dx \to 0 \)
\begin{align*}
A_1 \to \sum_{k\in\Ind} \int_0^{\infty} \int_{D^k} \big|u^k - c^k\big| \phi_t^k \,dx \,dt,
\end{align*}
and for \( \Dx \to 0 \)
\begin{align*}
A_2 \to \sum_{k\in\Ind} \big|u_0^k - c^k\big|(x) \phi^k(x, 0) \,dx, \qquad
A_4 \to 0.
\end{align*}
Thus, we are left with $A_3$. Since the scheme is the upwind method, we can write
\begin{align*}
A_3 &= \Dt \Dx \sum_{n=0}^{M} \biggl(\begin{aligned}[t]&\sum_{k\in\Indin} \sum_{i \in D_\disc^k} q_{c^k}^k \bigl( u_i^{k,n} \bigr) \biggl( \frac{\phi_{i+1}^{k,n} - \phi_{i}^{k,n}}{\Dx} \biggr) \\
&+ \sum_{k\in\Indout} \sum_{i \in D_\disc^k} q_{c^k}^k \bigl( u_{i}^{k,n} \bigr) \biggl( \frac{\phi_{i}^{k,n} - \phi_{i-1}^{k,n}}{\Dx} \biggr)\biggr)
\end{aligned}\\
&\to \sum_{k\in\Ind} \int_0^T \int_{D^k} q_{c^k}^k \bigl( u^k(x,t) \bigr) \phi_x (x,t) \,dx \,dt
\end{align*}
as \( \Dt, \Dx \to 0 \), due to the a.e.~pointwise convergence of $\bu_{\Dt}$ to $\bu$.
\end{proof}

To show compactness we want to apply Helly's theorem:
\begin{theorem}[Helly's theorem, {\cite[p.~437]{HoldenRisebro2015}}]
Let \( A\subset C([a,b]) \) be a set of functions on an interval $[a,b]\subset\R$ for which there exists some $M>0$ such that
\[
\TV(v) + \|v\|_\infty \leq M \qquad \forall\ v\in A.
\]
Then $A$ is relatively compact with respect to uniform convergence.
\end{theorem}

Now we have everything in place to proof a compactness theorem.

\begin{theorem}[Compactness and Convergence to a Weak Solution]
Fix \( T > 0 \). Assume that each flux function \( f^k \) is locally Lipschitz continuous and strictly monotone. Let $\germ_\disc^0$ be a set of discrete stationary states for the upwind method. Let \( \bu_\Dt \) be computed from the upwind method with initial data $\bar \bu \in \initialdataset(\germ_\disc^0) \cap L^1(\Omega;\lambda)$, and assume that $\TV(\bar \bu)<\infty$. Then the numerical solution \( \{\bu_\Dt\}_{\Dt>0} \) converges in $C([0,T],L_\loc^1(\Omega;\lambda))$ to a weak solution $\bu$.
\end{theorem}

\begin{proof}
We first show convergence of the sequence of functions \( \bg_\Dt\colon\Omega\times[0,T]\to\R \),
\[
\bg_\Dt(x,k,t)\coloneqq f^k(u^k_\Dt(x,t)).
\]
The sequence $\bg_\Dt$ is uniformly $L^\infty$ bounded, by Proposition~\ref{prop:linfbound}, and it is Lipschitz continuous in time:
\begin{align*}
\int_{\Omega} \big|\bg_\Dt(t^{n+1})-\bg_\Dt(t^n)\big|\,d\lambda & \leq C_f\int_{\Omega}\big|\bu_\Dt(t^{n+1})-\bu_\Dt(t^n)\big|\,d\lambda \\
& \leq C_f(C\TV(\bar\bu)+\bar M)\Dt,
\end{align*}
by Corollary~\ref{cor:timeContinuity}. We can bound the total variation of $\bg_\Dt$ by
\begin{align*}
\TV(\bg_\Dt(\cdot,t)) &= \sum_{k\in\Ind}\sum_{i \in D_\disc^k} \big|f^k( u_i^{k,n}) - f^k( u_{i-1}^{k,n})\big| \\
&\quad + \sum_{k\in\Indin}\big|f^k(u_0^n)-f^k(u_{-1}^{k,n})\big| \\
&\leq \sum_{k\in\Ind}\sum_{i \in D_\disc^k} \big|f^k( u_i^{k,n}) - f^k( u_{i-1}^{k,n})\big| + N\|\bg_\Dt\|_{\infty} \\
&= \sum_{k\in\Ind} \sum_{i \in D_\disc^k} \big|F_\iphf^{k,n} - F_\imhf^{k,n}\big| + N\|\bg_\Dt\|_{\infty} \\
&= \frac{\Dx}{\Dt}\sum_{k\in\Ind} \sum_{i \in D_\disc^k}\big|u_i^{k,n+1}-u_i^{k,n}\big| + N\|\bg_\Dt\|_{\infty} \\
&\leq \Dx(C\TV(\bar\bu) + \bar M) + N\|\bg_\Dt\|_{\infty}.
\end{align*}
Applying Ascoli's compactness theorem together with Helly's theorem,
we get the existence of a subsequence $\Dt_\ell\to 0$ such that $\bg_{\Dt_\ell} \to \bg$ in $C([0,T],L^1_\loc(\Omega;\lambda))$ for some function $\bg$. The strict monotonicity of $f^k$ implies that
\[
\bu_\Dt(x,k,t) = ( f^k )^{-1}\big(\bg_\Dt(x,k,t)\big),
\]
and hence, also $\bu_{\Dt_\ell}$ converges in $C([0,T],L^1_\loc(\Omega;\lambda))$ to some function $\bu$. Theorem~\ref{thm:laxwendroff} implies that $\bu$ is the entropy solution; since this solution is unique (Theorem~\ref{thm:entrsolnunique}), the entire sequence $\{\bu_\Dt\}_{\Dt>0}$ must converge to $\bu$.
\end{proof}

\section{Discrete Stationary Solutions for Monotone Flux Functions}\label{Sec:StationarySol}
So far we have shown that if a sufficiently large class of stationary
and discrete stationary solutions exists, then our equations on the
network are well posed and the finite volume numerical approximations
converge to the  entropy solution. In this section we show that such classes exist in the case where either all fluxes \( f^k \) are strictly increasing or all are strictly decreasing. We henceforth assume that all fluxes are increasing; one can attain analogous results for decreasing fluxes following the same arguments. In the following we want to investigate the sets of discrete stationary solutions implied by the upwind method.

We define
\begin{equation*}
  f_\into(u) \coloneqq \sum_{k\in\Indin} f^k(u), \qquad f_\out(u) \coloneqq
  \sum_{k\in\Indout} f^k(u) \qquad \text{for } u\in\R.
\end{equation*}
It is clear that \( f_\into, f_\out \) are monotone by the monotonicity of their summand components. In particular, the two functions are invertible.

For the upwind method the conditions \eqref{eq:consistencyInEdge} and \eqref{eq:consistencyOutEdge} become
\begin{subequations}
	\begin{align}
	f^k ( c^k ) &= f^k ( c^k ) \qquad \text{for } k\in\Indin,  \label{eq:consistencyInEdgeUpwind} \\
	f^k ( c^0 ) &= f^k ( c^k ) \qquad \text{for }k\in\Indout. \label{eq:consistencyOutEdgeUpwind}
	\end{align}
\end{subequations}
This is equivalent to
\begin{subequations}\label{eq:monotonedss}
	\begin{align}
	c^k &= c^k \qquad \text{for } k\in\Indin, \\
	c^0 &= c^k \qquad \text{for }k\in\Indout,
	\end{align}
\end{subequations}
due to the invertibility of the flux functions \( f^k \).
It is obvious as well, that for two different discrete stationary solutions \( \bc, \bd \) satisfying \( c^k \leq d^k \) for \( k \in \Ind \), we also have $c^0 \leq d^0$. Henceforth, we denote
\[
\germ_\disc^0 \coloneqq \bigl\{\text{all discrete stationary states for the upwind method}\bigr\}
\]
and we let
\[
\germ_\disc \coloneqq \bigl\{\big(c^{-N_\into},\dots,c^{-1},c^1,\dots,c^{N_\out}\big) \mid \bc\in\germ_\disc^0)\bigr\}.
\]
Although it might be too difficult to find a full characterization of the set $\initialdataset(\germ)$ of admissible initial data, we will be able to characterize large subsets of $\initialdataset(\germ)$. Let
\[
I_\into \coloneqq f_\into^{-1} (R_\into\cap R_\out), \qquad
I_\out \coloneqq f_\out^{-1} (R_\into\cap R_\out)
\]
where
\[
R_\into \coloneqq f_\into(\R), \qquad R_\out \coloneqq f_\out(\R).
\]
By the continuity of $f_\into, f_\out$, the sets $I_\into, I_\out$ are closed intervals.

\begin{theorem}
We have $\mathcal{L}\subset\initialdataset(\germ_\disc)$, where
\[
\mathcal{L}\coloneqq \Bigl\{ \bu\in L^\infty(\Omega;\lambda) \mid u^k(x) \in I_\into\ \forall\ k\in\Indin,\ u^k(x)\in I_\out\ \forall\ k\in\Indout \Bigr\}.
\]
In particular, if $f_\into, f_\out$ have the same range $R_\into,R_\out$, then $\initialdataset(\germ_\disc)=L^\infty(\Omega;\lambda)$.
\end{theorem}
\begin{proof}
Let $\bu\in\mathcal{L}$. Since $\bu\in L^\infty(\Omega;\lambda)$, and $I_\into,I_\out$ are closed, we also have
\[
\underline{c}_\into\coloneqq\inf_{\substack{x\in D^k\\ k\in\Indin}}u^k(x) \in I_\into, \qquad
\overline{c}_\into\coloneqq\sup_{\substack{x\in D^k\\ k\in\Indin}}u^k(x) \in I_\into
\]
and likewise for $\underline{c}_\out, \overline{c}_\out$. By continuity of $f_\into,f_\out$, there are $\underline{d}_\into\in I_\into$ and $\underline{d}_\out\in I_\out$ satisfying $\underline{d}_\into\leq\underline{c}_\into$ and $\underline{d}_\out\leq \underline{c}_\out$ so that $f_\into(\underline{d}_\into)=f_\out(\underline{d}_\out)$, that is, the vector $\underline{\bd}\coloneqq\bigl(\underline{d}_\into, \dots, \underline{d}_\into, \underline{d}_\out, \dots, \underline{d}_\out\bigr)$ is a stationary solution. This stationary solution clearly satisfies \eqref{eq:monotonedss}, whence $\bd\in\germ_\disc^0$.

In a similar way one finds a stationary solution $\overline{\bd}\coloneqq\bigl(\overline{d}_\into, \dots, \overline{d}_\into, \overline{d}_\out, \dots, \overline{d}_\out\bigr) \in \germ_\disc^0$ which bounds $\bu$ from above. Since now
\[
\underline{d}_\into\leq u^k(x)\leq\overline{d}_\into \;\forall\ k\in\Indin \qquad \text{and} \qquad
\underline{d}_\out\leq u^k(x)\leq\overline{d}_\out \;\forall\ k\in\Indout
\]
we conclude that $\bu\in \initialdataset(\germ_\disc)$.
\end{proof}

\begin{proposition}\label{prop:DiscGermUpwind}
Consider a conservation law on a network with strictly increasing fluxes \( f^k \). Let \( \germ_\disc^0 \) denote the set of all discrete stationary solutions for the upwind method. Then the set
\[ \germ_\disc \coloneqq \Big\lbrace \big(c^{-\Nin}, \dots, c^{-1}, c^1, \dots, c^{\Nout}\big) \mid \bc \in \germ^0_\monotone \Big\rbrace
\]
is a mutually consistent, maximal set of stationary solutions.
\end{proposition}
\begin{proof}
Every $\bc\in\germ_\disc$ is a stationary solution due to \eqref{eq:RHdss}.

To prove mutual consistency of \( \germ_\disc \) we plug a discrete stationary solution \( \bd\in\germ_\disc^0 \) into \eqref{eq:DiscEntropyIneqEdge} to get for \( n \in \N \),
\[
Q_{-\hf}^{k,n} \geq Q_{-\thf}^{k,n} \;\text{for }k \in \Ind_\into \quad\text{and}\quad
Q_{\thf}^{k,n} \geq Q_\hf^{k,n} \;\text{for } k \in \Ind_\out.
\]
Since we are using the upwind scheme, this reduces to
\[
q_{c^k}^k (d^k) \geq q_{c^0}^k (d^0) \;\text{for } k \in \Ind_\into \quad\text{and}\quad q_{c^0}^k (d^0) \geq q^k_{c^k} (d^k) \text{for } k \in \Ind_\out.
\]
In the same manner, plugging \( d^0 \) into \eqref{eq:DiscEntropyIneqVertex} gives us
\[ \sum_{k\in\Indin} Q_{-\hf}^{k,n} \geq \sum_{k\in\Indout} Q_\hf^{k,n}. \]
Combining these two observations, we get
\[ \sum_{k\in\Indin} q_{c^k}^k (d^k) \geq \sum_{k\in\Indin} q_{c^0}^k (d^0) \geq \sum_{k\in\Indout} q_{c^0}^k (d^0) \geq \sum_{k\in\Indout} q_{c^k}^k (d^k). \]
As $\bc,\bd$ were arbitrary, it follows that $\germ_\disc$ is mutually consistent.

If for some vector \( \bd \), \(\germ_\disc \cup \lbrace \bd \rbrace \) is mutually consistent, we know that
\[ \sum_{k\in\Ind_\into} q^k_{c^k}\big(d^k\big) \geq \sum_{k\in\Ind_\out} q^k_{c^k}\big(d^k\big). \]
We choose \( c^k = d^k\ \forall\ k \in \Ind_\into \) and \( c^0 = (f_\out)^{-1}(\sum_{k\in \Ind_\into} f^k(c^k)) \). Since all \( f^k \) are monotonically increasing, the entropy flux reduces to \( q^k_{c^k}(d^k)=\big|f^k(c^k) - f^k(d^k)\big| \) and thus,
\[ 0 \geq \sum_{k\in\Ind_\out} \big|f^k(c^0) - f^k(d^k)\big|, \]
which implies \( d^k = c^0 \) for \( k \in \Ind_\out \), and thus, $\bd\in\germ_\disc$. In other words, \( \germ_\disc \) is maximal.
\end{proof}

\begin{remark}
Due to Proposition \ref{prop:DiscGermUpwind} we now know that the doubling of variables argument in Section \ref{Sec:Stability} can be applied in this case and thus, the stability result holds.
\end{remark}

\begin{example}
Consider the traffic flow Example~\ref{ex:traffic} with $f^k(v)=\alpha^kf(v/\alpha^k)$ for $\alpha^k>0$ and $f$ some fixed flux, say, the Burgers flux $f(v)=v^2/2$. Restrict the solutions to the monotone region $u^k\geq0$. Consider the Godunov method with the numerical flux \( F(a,b) = f(a) \). The Rankine--Hugoniot condition reads
\[\sum_{k\in\Indin} \frac{(c^k)^2}{\alpha^k} = \sum_{k\in\Indout} \frac{(c^k)^2}{\alpha^k}. \]
For the discrete stationary solutions, the additional conditions \eqref{eq:dssEdge} read
\begin{align*}
f^k(c^0) = f^k ( c^k ) \qquad \text{for } k\in\Indout
\end{align*}
(the condition for $k\in\Indin$ automatically holds). Since $f^k$ is strictly monotone on $\R_+$, we get \( c^k = c^0 \) for all $k\in\Indout$. Thus, the set of all discrete stationary solutions is
\[
\germ_\disc^0 = \Big\{\bc\in\R^{N+1}\ :\ c^k=c^0\ \forall\ k\in\Indout,\ {\textstyle \sum_{k\in\Indin} \frac{(c^k)^2}{\alpha^k} = \sum_{k\in\Indout} \frac{(c^0)^2}{\alpha^k}} \Big\}.
\]
\end{example}

\section{Numerical Examples}\label{Sec:Numerics}
We show numerical experiments for some example cases including results for linear and nonlinear as well as convex and concave fluxes. In all experiments we use a CFL number of $\hf$ -- that is, $\Dt$ is chosen so that there is equality in \eqref{eq:cfl}.  In all experiments we compute the experimental order of convergence (EOC) as \( p \approx \dfrac{\log(e^{j+1}/e^{j})}{\log(\Dx_{j+1}/\Dx_j)} \) on a series of successive grids with \( 2^{j} \) cells, where \( e^j \) denotes the \( L^1 \) error on grid level \( j \). The error is computed as the $L^1$ difference to a high-resolution reference solution. All errors and EOC are displayed in Table~\ref{tab:EoC}.

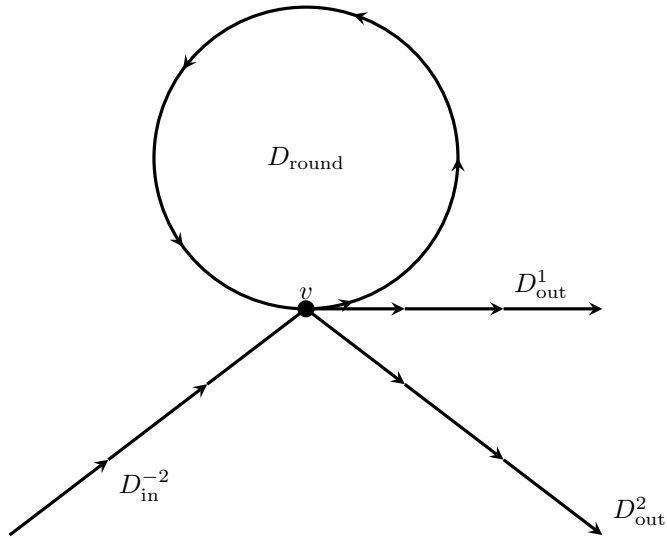
\begin{figure}[ht]
\centering
\begin{tikzpicture}[>=stealth]
\coordinate (origin) at (0.0,-3.5);
\coordinate (horizontalshift) at (0.1,0);
\coordinate (verticalshift) at (0,0.1);
\coordinate (offset1) at (0.5,0);
\coordinate (offset2) at (1.5,0);
    \draw[->, very thick] (origin) -- ($(origin) + 13*(horizontalshift)$);
    \draw[->, very thick] ($(origin) + 13*(horizontalshift)$) -- ($(origin) -- 26*(horizontalshift)$) node[anchor = south west] {$D_\out^1$};
    \draw[->, very thick] ($(origin) + 26*(horizontalshift)$) -- ($(origin) + 39*(horizontalshift)$);
    \draw[->, very thick] (origin) -- ($(origin) + 13*(horizontalshift) - 10*(verticalshift)$);
    \draw[->, very thick] ($(origin) + 13*(horizontalshift) - 10*(verticalshift)$) -- ($(origin) + 26*(horizontalshift) - 20*(verticalshift)$);
    \draw[->, very thick] ($(origin) + 26*(horizontalshift) - 20*(verticalshift)$) -- ($(origin) + 39*(horizontalshift) - 30*(verticalshift)$) node[anchor = south west] {$D_\out^2$};
	\draw[very thick,
		    decoration={markings, mark=at position 0.0 with {\arrow{>}}},
	        decoration={markings, mark=at position 0.2 with {\arrow{>}}},
	        decoration={markings, mark=at position 0.4 with {\arrow{>}}},
	        decoration={markings, mark=at position 0.6 with {\arrow{>}}},
	        decoration={markings, mark=at position 0.8 with {\arrow{>}}},
	        postaction={decorate}
	        ]
	        ($(origin)+20*(verticalshift)$) circle (2) node[anchor = center] {$D_{\textrm{round}}$};
    \draw[->, very thick] ($(origin) - 30*(verticalshift) - 39*(horizontalshift)$) -- ($(origin) - 20*(verticalshift) - 26*(horizontalshift)$) node[anchor = north west] {$D_\into^{-2}$};
    \draw[->, very thick] ($(origin) - 20*(verticalshift) - 26*(horizontalshift)$) -- ($(origin) - 10*(verticalshift) - 13*(horizontalshift)$);
    \draw[very thick] ($(origin) - 10*(verticalshift) - 13*(horizontalshift)$) -- (origin);
    \draw[fill=black] (origin) circle (.7ex) node[anchor = south] {$v$};

\end{tikzpicture}
\caption{A network with a periodic edge.}\label{fig:Roundaboutgraph}
\end{figure}
\begin{example}[Burgers' equation with roundabout]\label{ex:burgersroundabout}
In this example we include a roundabout -- an edge whose endpoints meet at the same vertex, as shown in Figure~\ref{fig:Roundaboutgraph}. We also include an ingoing edge and two outgoing edges, amounting to a total of two ingoing and three outgoing edges. As initial data we choose constants on the roundabout and the outgoing edges and two different constants on the independent ingoing edge. After a while the shock in the initial data on the independent ingoing edge will reach the edge and create new Riemann problems. We choose the initial data
\begin{align*}
\bar u^{-1}(x) = \begin{cases} 2 & \text{if } 0 \leq x < 0.5, \\ \sqrt{2} & \text{if } 0.5 \leq x, \end{cases}\qquad \bar u^{-2} = \bar u^1 = \bar u^2 = \bar u^3 \equiv 1.
\end{align*}
We take all edges to have length $1$ and choose zero Neumann boundary data on the outer boundaries. On the vertex we set \( u_0^0 = 1 \). On the ingoing edge with index \( -1 \) we have a travelling shock wave
\begin{equation*}
u^{-1}(x,t) = \begin{cases} 2 & \text{if } 0 \leq x < \frac{1}{2-\sqrt{2}} t \\
 \sqrt{2} & \text{if } \frac{1}{2-\sqrt{2}} t \leq x \end{cases}
\end{equation*}
which will hit the vertex at \( t^* = 1-\frac{1}{\sqrt{2}} \). To compute the solution after $t^*$ we compute the new vertex value \( c^0 = \sqrt{\nicefrac{5}{3}} \) and therefore get the Riemann problem
\begin{align*}
u^k(x,t^*) = \begin{cases} \sqrt{\nicefrac{5}{3}} & \text{if } x = 0, \\ 1 & \text{if } 0 < x \end{cases}
\end{align*}
for \( k = 1,2,3 \), which results in a travelling shock wave with speed \( s = \frac{1}{\sqrt{\nicefrac{5}{3}} - 1} \). At time \( t^{**} \coloneqq \sqrt{\nicefrac{5}{3}} - \frac{1}{\sqrt{2}} \) the travelling shock wave which originated on the roundabout edge hits the vertex once again, resulting in a new set of Riemann problems on the outgoing edge. This process will continue in a periodic fashion.

We compute up to time $T=0.5$. A plot of the exact and approximate solution to this example at two different times is shown in Figure~\ref{fig:RoundaboutBurgers}. The accuracy and order of convergence of the numerical approximation are shown in Table~\ref{tab:EoC}. 
\begin{figure}
 \includegraphics[scale=.42]{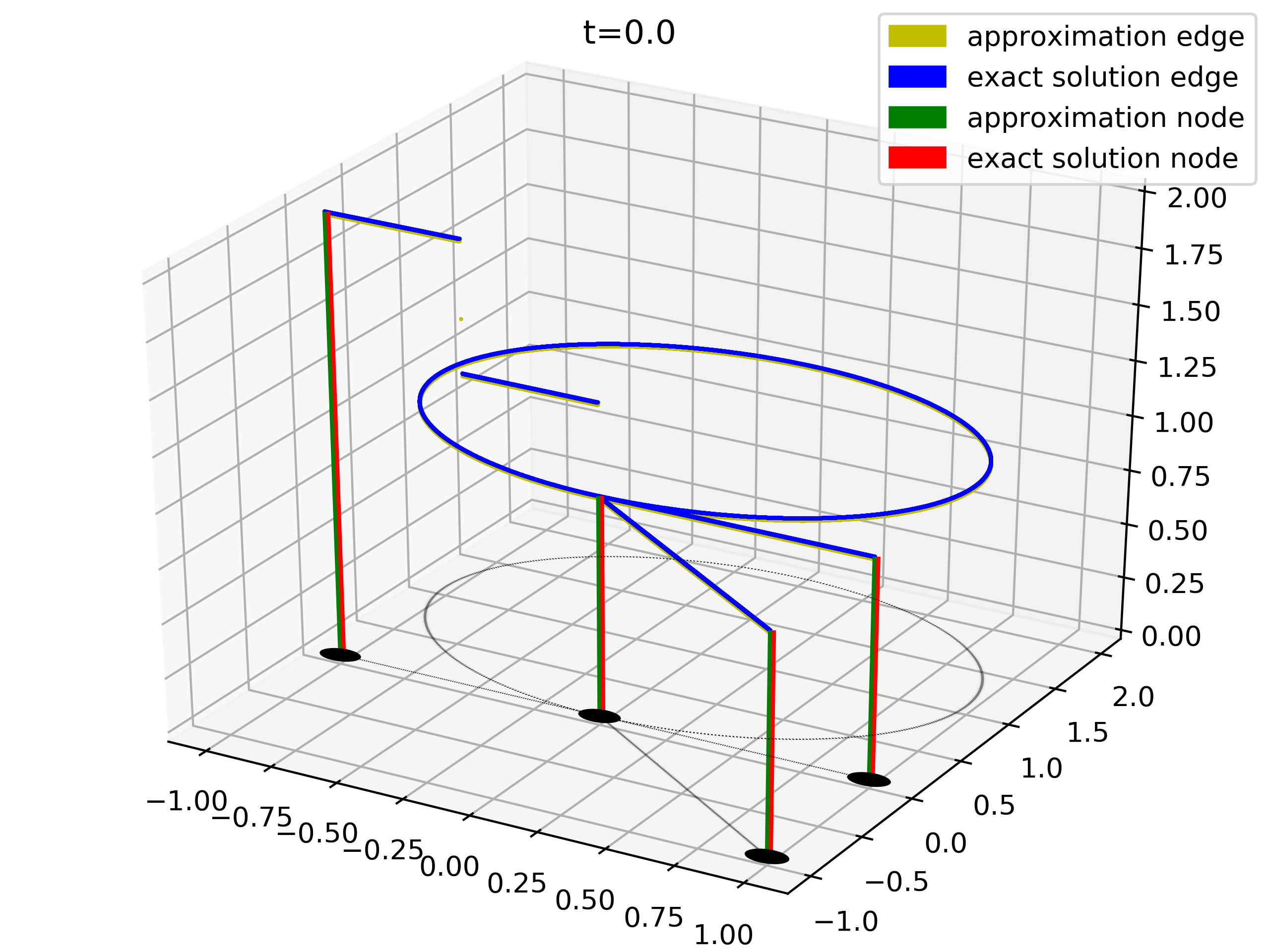}
 \includegraphics[scale=.42]{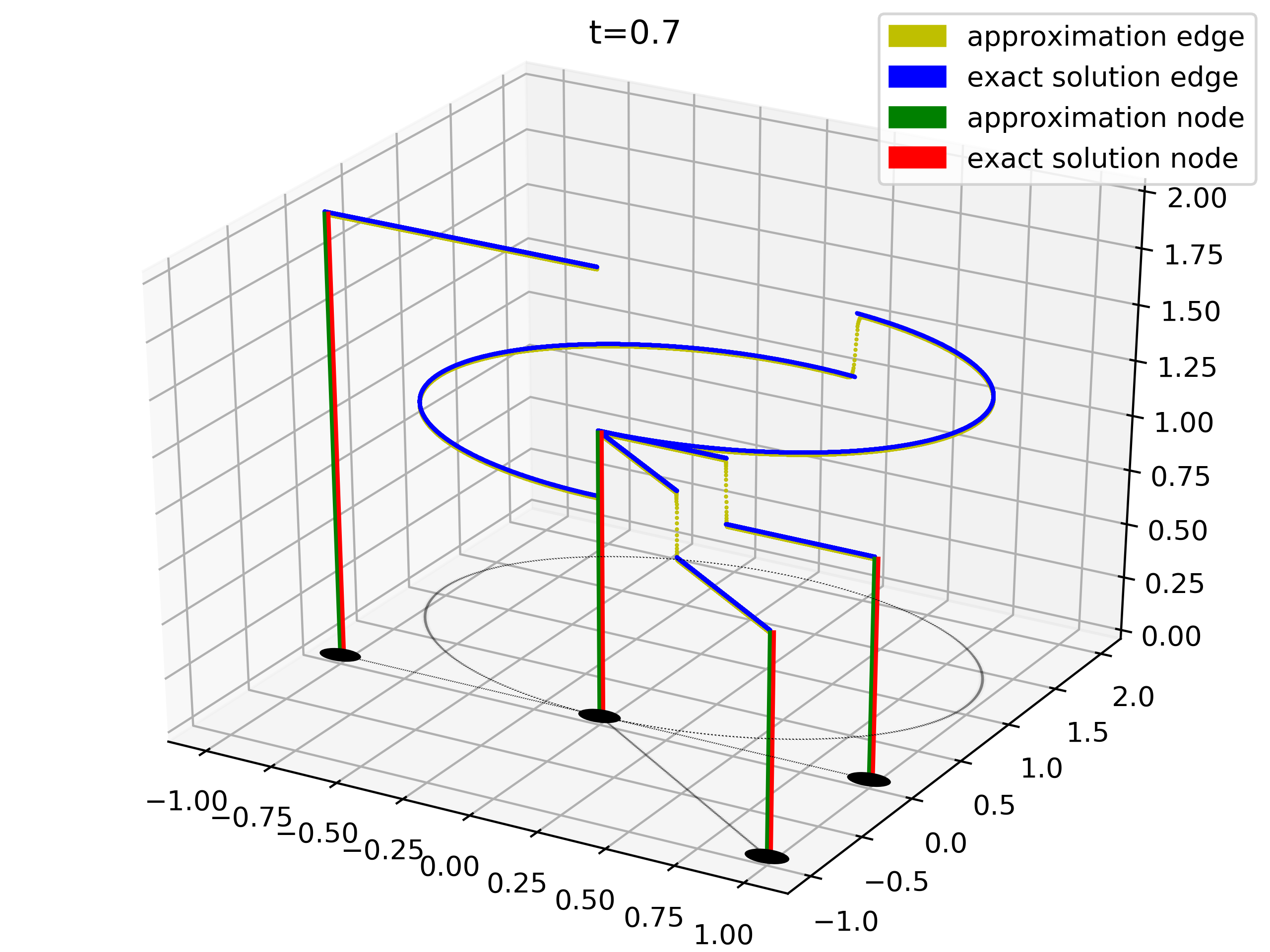}
 \caption{Initial state and state at $t=0.7$ of a {Burgers-type equation} with travelling shock wave which hits the vertex at time $t=1-\frac{1}{\sqrt{2}}$. Here, the graph includes a periodic edge.}\label{fig:RoundaboutBurgers}
\end{figure}
\end{example}

\begin{example}\label{ex:holdenrisebro}
We construct an example where we take the flux function from the traffic flow example in \cite{HoldenRisebro1995}, \( f(u) = 4u(1-u) \), but allow for different fluxes on different edges, \( f^k( u ) = \alpha^k f \big( \frac{u}{\alpha^k} \big) \) for $\alpha^k>0$, and compute on a star shaped graph with two ingoing edges and three outgoing edges like in Figure~\ref{fig:Stargraph}. The initial data is chosen so that all fluxes are strictly increasing over the range of $\bar u$; thus, the fluxes $f^k$ are in effect monotonously increasing functions. We choose constant solutions on the two ingoing roads and constant initial data on the outgoing roads which are chosen such that on one road a shock will develop, on one road the solution will stay constant over time and on one road a rarefaction wave will develop.

Solving the conditions \eqref{eq:RHdss}, \eqref{eq:dssEdge} for $c^0$ yields
\[ c^0 = \frac{3 \pm \sqrt{9 - 4 \bigl(\frac{1}{\alpha_1} + \frac{1}{\alpha_2} + \frac{1}{\alpha_3}\bigr) \big( u^{-1} - \frac{1}{\alpha_{-1}} ( u^{-1} )^2 + u^{-2} - \frac{1}{\alpha_{-2}} ( u^{-2} )^2 \big)}}{2\big(\frac{1}{\alpha_1} + \frac{1}{\alpha_2} + \frac{1}{\alpha_3}\big)}. \]
For the incoming edges to have a monotonically increasing flux we impose \( u_{-i} \leq \frac{1}{2} \alpha_{-i},\) for \( i=1,2 \) and for outgoing edges \( c^0 \leq \frac{1}{2} \min \lbrace \alpha_1, \alpha_2, \alpha_3 \rbrace \). We choose \( \alpha_{-1} = \alpha_{-2} = 1 \), \( \alpha_1 = \alpha_2 = 4 \) and \( \alpha_3 = 2 \) with initial data
\begin{align*}
\bar u^{-1} = \bar u^{-2} \equiv 0.5,\qquad \bar u^1 \equiv 0, \qquad \bar u^2 \equiv \frac{1}{2} \big(3-\sqrt{7}\big),\qquad \bar u^3 \equiv 1.
\end{align*}
This gives us \( c^0 = \frac{1}{2}(3 - \sqrt{7}) \). On the outer boundary we choose zero Neumann boundary conditions. For \( u^3 \) we will get a shock
\[ u^3(x,t) = \begin{cases} \frac{1}{2}(3 - \sqrt{7}) & \text{if } x < st,
\\ 0 & \text{if } x\geq st, \end{cases} \]
with speed \(s= \frac{(3-\sqrt{7})(2-\frac{3-\sqrt{7}}{4})-3}{\frac{3-\sqrt{7}}{2}-1},\) and a rarefaction wave for \( u^1 \) of the form
\[
u^1(x,t) = \begin{cases} \frac{3-\sqrt{7}}{2} &  \text{if }x < 2 (\sqrt{7}-1)t, \\
1-\frac{x}{4t} & \text{if } 2(\sqrt{7}-1)t \leq x < 4t, \\
0 & \text{if } 4t \leq x. \end{cases}
\]
On edge $2$ we get the constant solution \( u^2 \equiv \frac{1}{2}(3 - \sqrt{7}) \).

We compute up to time $T=0.2$. A plot of the exact and approximate solution at two different timepoints is shown in Figure~\ref{fig:MonotonousTrafficFlow}. Accuracy and order of convergence of the numerical approximation are shown in Table~\ref{tab:EoC}.
\begin{figure}
 \includegraphics[scale=.42]{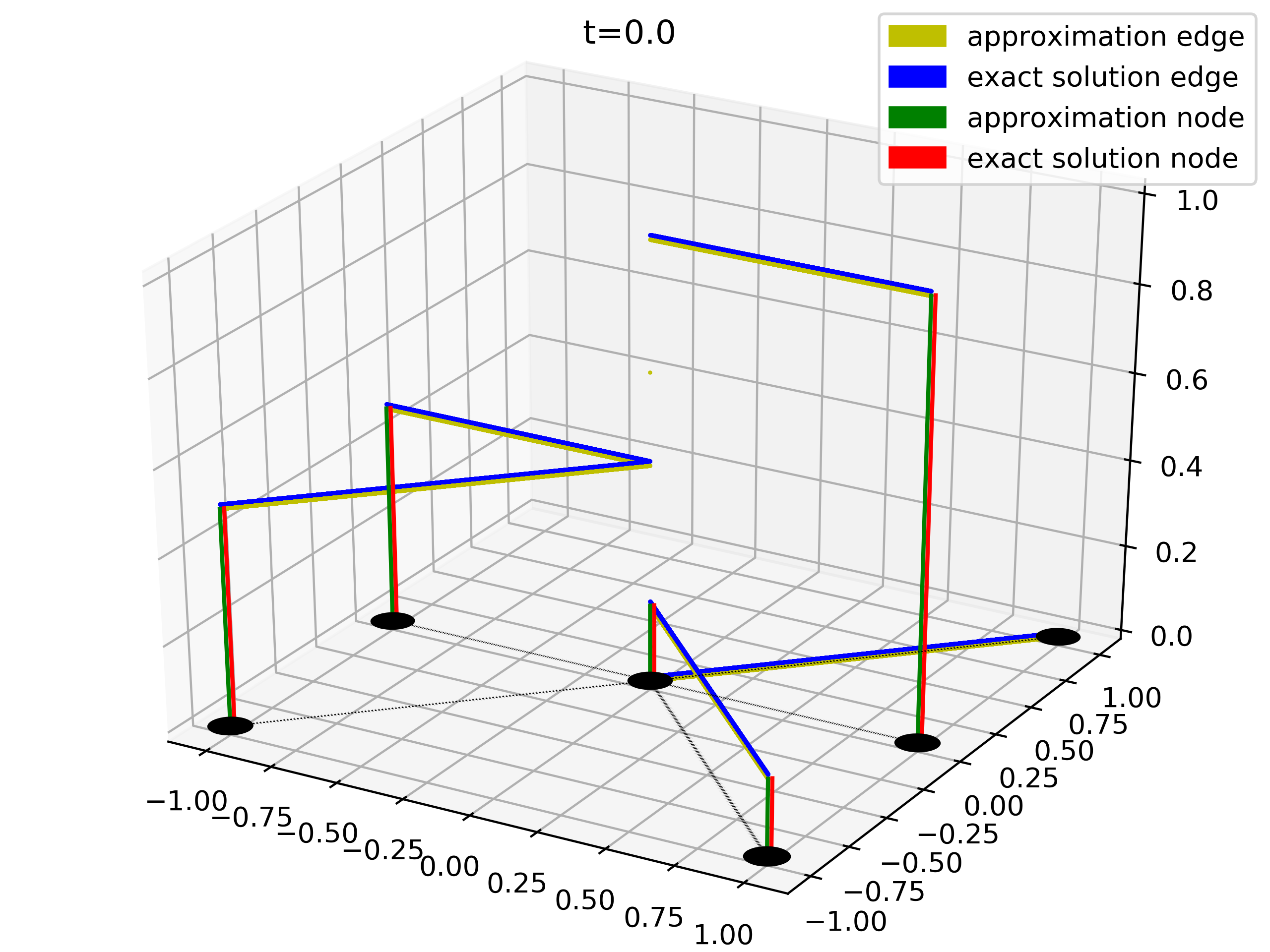}
 \includegraphics[scale=.42]{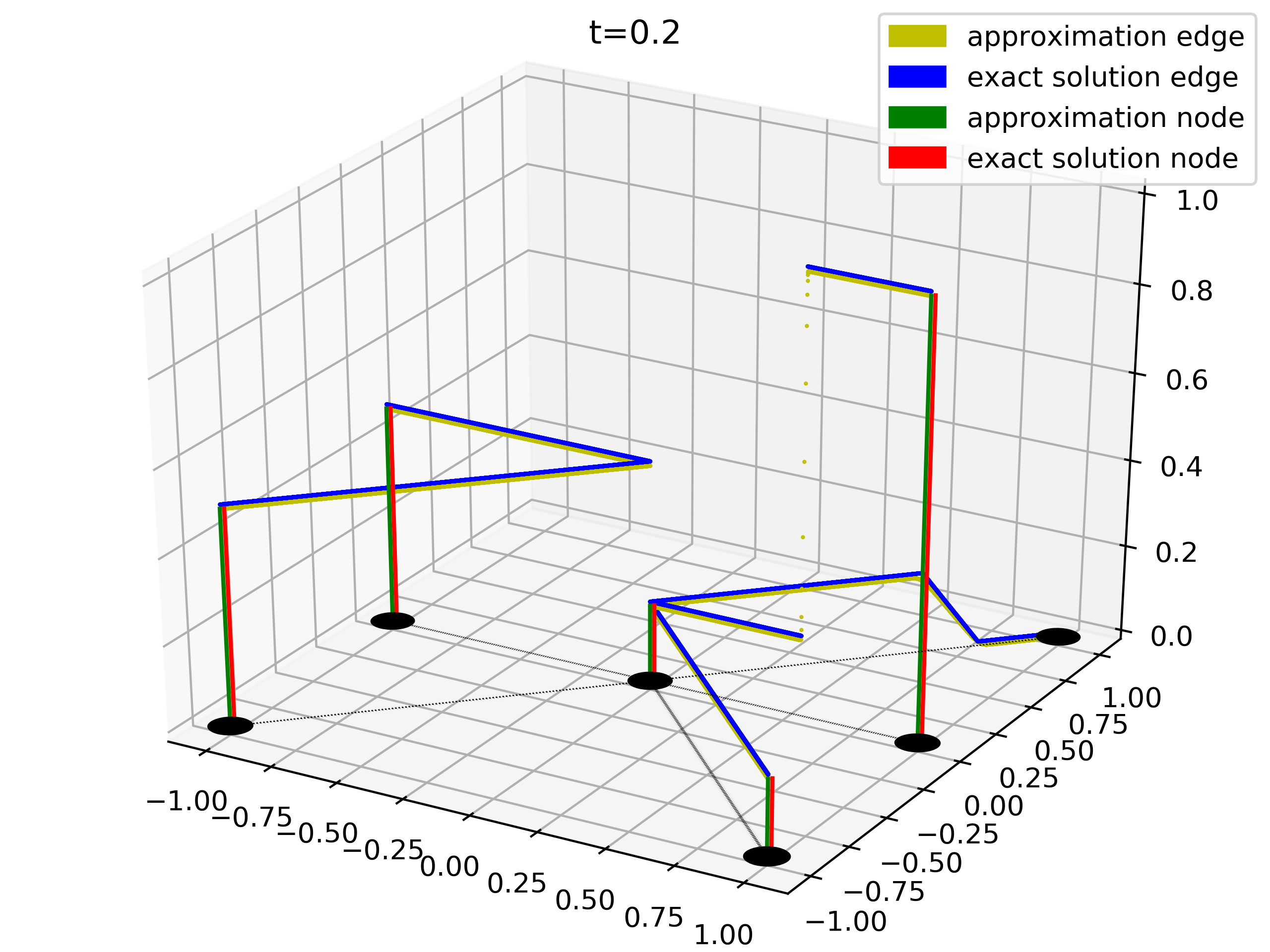}
 \caption{Initial state at $t=0$ and state at $t=0.2$ of a \textit{traffic flow problem} with an initial shock at the vertex developing a different elementary wave on each outgoing edge.}\label{fig:MonotonousTrafficFlow}
\end{figure}
\end{example}

In addition to the examples described above we show errors and experimental order of convergence (EOC) for several additional examples in Table~\ref{tab:EoC}.

\begin{example}[EOC: Linear advection]\label{ex:linadv}
We consider a linear advection equation with two ingoing edges and three outgoing edges as in Figure~\ref{fig:Stargraph} with initial data
\begin{align*}
\bar u^{-1}(x,t) = \begin{cases} 2 & 0 \leq x < 0.8, \\ 1 & 0.8 \leq x, \end{cases}\qquad \bar u^{-2} \equiv 1,\qquad \bar u^1 = \bar u^2 = \bar u^3 \equiv \frac{2}{3},
\end{align*}
and Dirichlet boundary conditions adapted to the edge values. We initialize the vertex node by $u_0^0 = \frac{2}{3}$. We compute up to time $T=0.5$.
\end{example}

\begin{example}[EOC: Burgers' equation with elementary waves]\label{ex:burgerselementary}
We choose \( \bar u^{-1} = \bar u^{-2} \equiv 1 \) as initial data on the ingoing roads and \( \bar u^1 \equiv 0 \), \( \bar u^2 \equiv\sqrt{\nicefrac{2}{3}}\) and \( \bar u^3 \equiv 2 \) on the outgoing edges of a star shaped graph as in Figure~\ref{fig:Stargraph}. The conditions on the numerical flux imply then \( 3(c^0)^2 = 2  \Leftrightarrow c^0 = \sqrt{\nicefrac{3}{2}} \). Thus, we get the following Riemann problems on the outgoing roads:
\begin{align*}
\bar u^1(x) = \begin{cases} \sqrt{\nicefrac{2}{3}} & x = 0, \\ 0 & x > 0, \end{cases}\qquad \bar u^2 = \begin{cases} \sqrt{\nicefrac{2}{3}} & x = 0, \\ \sqrt{\nicefrac{2}{3}} & x > 0, \end{cases}\qquad \bar u^3 = \begin{cases} \sqrt{\nicefrac{2}{3}} & x = 0, \\ 2 & x > 0, \end{cases}
\end{align*}
with zero Neumann boundary conditions at the outer edges. The solution to these problems are a shock, a constant solution and a rarefaction wave, respectively. We compute up to time $T=0.3$.
\end{example}

\begin{example}[EOC: Burgers' equation with travelling shock]\label{ex:burgersshock} We consider a Burgers-type equation with two ingoing edges and three outgoing edges as in Figure~\ref{fig:Stargraph} with initial data
\begin{align*}
\bar u^{-1}(x) = \begin{cases} 2 & 0 \leq x < 0.8, \\ 1 & 0.8 \leq x, \end{cases}\qquad \bar u^{-2} \equiv 1,\qquad \bar u^1 = \bar u^2 = \bar u^3 \equiv \sqrt{\frac{2}{3}},
\end{align*}
with Dirichlet boundary conditions of the same value as the associated edge. On the vertex node the initial condition is chosen as $u_0^0 = \sqrt{\frac{2}{3}}$. We compute up to $T=0.5$.
\end{example}

\begin{table}[h]
\centering
\setlength{\tabcolsep}{3pt}
\renewcommand{\arraystretch}{1.3}
\scriptsize
\begin{tabular}{c|c|c|c|c|c|c|c|c|c|c}
\hline
 & \multicolumn{2}{c|}{Example~\ref{ex:linadv}} & \multicolumn{2}{c|}{Example~\ref{ex:burgersshock}} & \multicolumn{2}{c|}{Example~\ref{ex:burgerselementary}}  & \multicolumn{2}{c|}{Example~\ref{ex:burgersroundabout}} & \multicolumn{2}{c}{Example~\ref{ex:holdenrisebro}} \\
\hline
Grid level & $L^1$ error   & EOC  & $L^1$ error & EOC & $L^1$ error & EOC & $L^1$ error & EOC & $L^1$ error & EOC \\
\hline
 3 & 0.10877 &   -  & 0.11630 &   -  & 0.14459 & - & 0.07087 & - & 0.09904 & - \\
 4 & 0.05496 & 0.98 & 0.07136 & 0.70 & 0.08016 & 0.85 & 0.0546  & 0.38 & 0.04913 & 1.01 \\
 5 & 0.03649 & 0.59 & 0.04372 & 0.71 & 0.04651 & 0.79 & 0.03117 & 0.81 & 0.02844 & 0.79 \\
 6 & 0.02629 & 0.47 & 0.02255 & 0.96 & 0.02711 & 0.78 & 0.01903 & 0.71 & 0.01627 & 0.81 \\
 7 & 0.01830 & 0.52 & 0.01360 & 0.73 & 0.01495 & 0.86 & 0.01115 & 0.77 & 0.00919 & 0.82 \\
 8 & 0.01255 & 0.54 & 0.00653 & 1.06 & 0.00925 & 0.69 & 0.00644 & 0.79 & 0.00527 & 0.80 \\
 9 & 0.00883 & 0.51 & 0.00325 & 1.01 & 0.00480 & 0.95 & 0.00330 & 0.96 & 0.00268 & 0.98 \\
10 & 0.00625 & 0.50 & 0.00160 & 1.02 & 0.00295 & 0.70 & 0.00173 & 0.93 & 0.00150 & 0.84 \\
11 & 0.00442 & 0.50 & 0.00086 & 0.90 & 0.00152 & 0.96 & 0.00085 & 1.03 & 0.00084 & 0.84 \\
12 & 0.00312 & 0.50 & 0.00040 & 1.10 & 0.00081 & 0.91 & 0.00042 & 1.02 & 0.00047 & 0.84 \\
\hline
\end{tabular}
\caption{$L^1$ errors and estimated orders of convergence (EOC) for a selection of examples.}\label{tab:EoC}
\end{table}

\subsection{Comments on the experiments}
Convergence order estimates for finite volume methods for nonlinear scalar conservation laws are due to Kuznetsov \cite{Kuznetsov76} for the continuous flux case and due to Badwaik, Ruf \cite{BadwaikRuf2020} for the case of monotone fluxes with points of discontinuity. In both of those cases the analytically proven convergence rate is at least \( \sqrt{\Dx} \). Our numerical experiments indicate the same lower bound on the convergence rate for our numerical methods on graphs. Considering the fact that \( f_\into \) and \( f_\out \) from Section~\ref{Sec:StationarySol} are monotone it might be possible to generalize the result of Badwaik and Ruf to networks.

\section{Summary and Outlook}
In conclusion we have defined a framework for the analysis and numerical approximation of conservation laws on networks.
We extended the concepts well known from the conventional case such as weak solution, entropy solution and monotone methods to make sense on a directed graph.
We defined a reasonable entropy condition under which we have shown stability and uniqueness of an analytic solution.
Existence is shown by convergence of a conservative, consistent, monotone difference scheme.
In an upcoming work \cite{FjordholmMuschRisebro2020} we want to address convergence of a numerical method where the fluxes $f^k$ are not monotone but concave, as is usually found in traffic flow models. This includes deriving a sufficiently large set of stationary and discrete stationary solutions for this case. Further, we want to extend our model to include boundary conditions and derive a convergence order estimate for numerical approximations.
As for future work, a generalization to systems of conservation laws would be highly desirable. One could also try to construct numerical schemes for equations incorporating diffusive fluxes like it was done in \cite{KarlsenRisebroTowers2002} on the line.
Generalized models would span more complex scenarios such as blood circulation \cite{BressanCanicGaravelloPiccoli2013} in a network of veins or a river delta by the means of Euler equations and shallow water equations, respectively.
\bibliographystyle{plain}
\bibliography{bibliography.bib}
\end{document}